\documentclass[11pt]{article}
\usepackage{amsmath}
\usepackage{amssymb,amsmath}
\usepackage{mathrsfs}
\usepackage{graphicx}
\let\r=\rho

\let\f=\frac
\let\p=\psi

\let\D=\Delta

\let\wt=\widetilde

\def\cA{{\cal A}}

\def\cC{{\cal C}}

\def\cF{{\cal F}}

\def\cP{{\cal P}}

\def\cS{{\cal S}}

\def\cZ{{\cal Z}}

\def\na{\nabla}

\def\p{\partial}

\def\dv{\mbox{div}}

\def\eqdefa{\buildrel\hbox{\footnotesize def}\over =}

\def\C{\mathop{\bf C\kern 0pt}\nolimits}
\def\DD{\mathop{\bf D\kern 0pt}\nolimits}
\def\K{\mathop{\bf K\kern 0pt}\nolimits}
\def\N{\mathop{\bf N\kern 0pt}\nolimits}
\def\Q{\mathop{\bf Q\kern 0pt}\nolimits}
\def\R{\mathop{\bf R\kern 0pt}\nolimits}

\def\ddq{\dot \Delta_q}

\renewcommand{\div}{\mbox{\rm div}\;\!}

\newcommand{\Dv}{{\rm div}}

\newcommand{\beq}{\begin{equation}}
\newcommand{\eeq}{\end{equation}}
\newcommand{\ben}{\begin{eqnarray}}
\newcommand{\een}{\end{eqnarray}}
\newcommand{\beno}{\begin{eqnarray*}}
\newcommand{\eeno}{\end{eqnarray*}}

\newtheorem{Theorem}{Theorem}[section]
\newtheorem{Definition}[Theorem]{Definition}
\newtheorem{Proposition}[Theorem]{Proposition}
\newtheorem{Lemma}[Theorem]{Lemma}

\newtheorem{Remark}[Theorem]{Remark}

\numberwithin{equation}{section}

\allowdisplaybreaks \numberwithin{equation} {section}

\begin{document}
\title{Existence, uniqueness and optimal decay rates for the 3D   compressible Hall-magnetohydrodynamic system
     \thanks {Research supported by the
National Natural Science Foundation of China (11501332,11171034,11371221), the  Natural Science Foundation of Shandong Province (ZR2015AL007),
 and Young Scholars Research Fund of Shandong University of Technology.}
}
\author{ Fuyi  Xu$^{a\dag}$ \ \
         Meiling  Chi$^{a}$ \ \ Lishan   Liu$^{b,c}$ \ \ Yonghong  Wu$^{c}$  \\[2mm]
 { \small $ ^a$ School of Mathematics and  Statistics, Shandong University of Technology,}\\
  { \small Zibo,    255049,  Shandong,    China}\\
{\small $ ^b$ School  of Mathematics Science, Qufu normal
University,} \\ { \small Qufu,  263516, Shandong, China}\\
{\small $ ^c$ Department of Mathematics and Statistics, Curtin
University,}\\
{\small Perth, 6845, WA,  Australia}}
         \date{}
         \maketitle
\noindent{\bf Abstract}\ \ \ We are concerned with the study of
the Cauchy problem to  the 3D  compressible Hall-magnetohydrodynamic system. We first establish the  unique global solvability  of  strong
solutions to the system when  the initial data are  close to a stable
equilibrium state in critical Besov spaces. Furthermore,  under a suitable additional condition involving only the low frequencies of the  data and in $L^{2}$-critical regularity framework,  we  exhibit the optimal time decay rates  for  the  constructed global
solutions. The proof relies on an application of  Fourier analysis to a mixed  parabolic-hyperbolic system, and on a refined time-weighted energy functional.
\vskip   0.2cm \footnotetext[1]{$^\dag$Corresponding author.}
\vskip   0.2cm \footnotetext[2]{E-mail addresses: zbxufuyi@163.com(F. Xu),\ chimeiling0@163.com(M. Chi),\ mathlls@163.com(L. Liu),\  y.wu@curtin.edu.au(Y. Wu).} \setlength{\baselineskip}{20pt}

\section{Introduction and main results}
\setcounter{section}{1}\setcounter{equation}{0} \ \  In this paper, we consider the following
3D compressible Hall-magnetohydrodynamic equations \cite{ADFL}:
\begin{align}\label{equ:CMHD}
\left\{
\begin{aligned}
&\p_t\rho+\textrm{div}(\rho u)=0,\\
&\p_t(\rho u)+\textrm{div}(\rho u\otimes u)-\mu\Delta u-(\lambda+\mu)\na\textrm{div}u+\na p=(\nabla\times H)\times H, \\
&\p_tH-\nabla\times(u\times H)+\nu \nabla\times(\nabla\times H)+\nabla\times \Big(\frac{(\nabla\times H)\times H}{\rho}\Big)=0, \\
&\textrm{div}H=0,
\end{aligned}
\right.
\end{align}
where $\rho(t,x)$, $u(t,x)$, $H(t,x)$ denote, respectively, the
density, velocity,  and magnetic field. $p=p(\rho)$ is  pressure
satisfying $p'(\rho)>0$ and for all $\rho>0$. The Lam\'e
coefficients $\mu$ and $\lambda$  satisfy the physical conditions
\begin{equation}\label{condition}\mu>0,\ \ 2\mu+3\lambda> 0,\end{equation}
which ensures that the operator $-\mu\Delta-(\lambda+\mu)\nabla \Dv$
is a strongly elliptic operator and $\nu>0$ is the magnetic diffusivity acting as a
magnetic diffusion coefficient of the magnetic field. Here, we simply set
$\nu=1$ since its size does not play any role in our analysis. In this
paper, we are concerned with the Cauchy problem of the system
\eqref{equ:CMHD} in $\mathbb{R}_{+}\times \mathbb{R}^3$ subject to the
initial data
\begin{equation}\label{initial value} (\rho,u,H)|_{t=0}=(\rho_0, u_0,
H_0).\end{equation}
In many current physics problems, Hall-MHD is required. For example, when magnetic shear is large, which precisely occurs during reconnection events, the influence of the Hall term becomes dominant. However, this term is usually small and can be neglected in laminar situations, which is why conventional MHD models ignore it.  When $\rho=const$,  system \eqref{equ:CMHD} becomes
incompressible Hall-MHD system, which has received many studies. The first systematic study of Hall-MHD is due to Lighthill \cite{MJ} followed by Campos \cite{LCC}. The Hall-MHD is indeed needed for such problems as magnetic reconnection in space plasmas \cite{HG}, star formation \cite{SB}, and neutron stars \cite{DU}. A physical review on these questions can be found in \cite{PM}. Mathematical derivations of Hall-MHD equations from either two-fluids or kinetic models can be found in \cite{ADFL} and in this paper, the first existence result of global weak solutions is given.  In recent years, a number of works have been dedicated to this  model, we refer readers to Refs \cite{ADFL,chae1, chae2, chae3,fan3,fan2,fan1,fan7,WZ1} for more discussions.  When $\Dv u=0$,  system (\ref{equ:CMHD}) becomes the
 density-dependent Hall-MHD system, which
has been investigated by many authors, and for more details, see \cite{fan4,fan6}.

When the Hall effect term $\nabla\times \Big(\frac{(\nabla\times H)\times H}{\rho}\Big)$
is neglected, system (\ref{equ:CMHD}) reduces to  well-known compressible MHD system.  There are many results regarding
the global existence of the solutions and the decay of the smooth solutions to
the compressible MHD equations, see \cite{SK,UKS,TW} and references therein.  However, to the best of our
knowledge, very few results have been established on the dynamics of
the global solutions to the 3D compressible Hall-MHD system, especially on the  temporal decay of the solutions.
Recently, Fan et al.  \cite{fan3} first
obtained the global existence and the optimal decay rates of  smooth solutions to the 3D compressible Hall-MHD equations
\eqref{equ:CMHD} where  the initial data are close to an
equilibrium state in $H^3(\mathbb{R}^{3}_{x})$ and belong to $
L^1(\mathbb{R}^{3}_{x})$. Later, Gao and Yao \cite{GY} improved the result from \cite{fan3}. They proved the global existence
of strong solutions by the standard energy method under the condition that the initial data are
close to the constant equilibrium state in the lower regular spaces $H^{2}(\mathbb{R}^{3}_{x})$ and obtained optimal decay rates for the constructed
global strong
solutions in $L^{2}$-norm  if the initial data belong to $L^{1}$ additionally.
More  recently,
Xu et al.  \cite{Xu3} proved the global existence and temporal  decay rates of  the
solutions to the system \eqref{equ:CMHD} when  the initial data are  close to a stable
equilibrium state in $H^3(\mathbb{R}^{3}_{x})\cap
\dot{H}^{-s}(\mathbb{R}^{3}_{x})$ for some $s\in \lbrack 0,3/2)$ by using a pure energy method.
Obviously, all these results  for  compressible Hall-MHD equations \eqref{equ:CMHD} need higher
smoothness for the initial data. The price to pay however is that
assuming higher smoothness precludes from using a critical function
spaces framework. To  our knowledge, few results have
been established for the 3D compressible Hall-MHD equations in critical spaces. In the present article, we  first study the global
well-posedness in critical regularity framework. In order to tackle the problem, we need make scaling analysis of the system \eqref{equ:CMHD}.
Different from  compressible MHD equations, Hall-term breaks the natural scaling of the Navier-Stokes equations.
In fact, if we set the
density  $\rho=const$ and the fluid velocity $u\equiv0$, then \eqref{equ:CMHD} is reduced to
\begin{equation}\label{Heq} \p_tH-\Delta H+\nabla\times \Big((\nabla\times H)\times H\Big)=0,
\quad \hbox{in}\quad \mathbb{R}_{+}\times \mathbb{R}^3.\end{equation}
Assume that magnetic field $H$ satisfy \eqref{Heq}, then the function $H_{l}(x,t)=H(lx,l^{2}t)$ for $l>0$ form a solution to \eqref{Heq} again.
So, $\nabla H$ of \eqref{Heq} has the same scaling with the fluid velocity $u$ to the  usual Navier-Stokes equations. Motivated by  the scaling of 3D compressible Navier-Stokes equations and \eqref{Heq}, we find that the function space
for the velocity $u$ is similar to $\dot{H}^{\f{1}{2}}(\mathbb{R}^{3}_{x})$ and the regularity of the density $c=\rho-\bar{\rho}$ and the magnetic field $H$ is one order higher than that of the velocity $u$,  and
hence a function space for $c$ and $H$ which is similar to $\dot{H}^{\f32}$ would be
a candidate. Unfortunately, the function space
$\dot{H}^{\f{3}{2}}(\mathbb{R}^{3}_{x})$ does not turn to be a good candidate for $c$ and $H$
since  $\dot{H}^{\f{3}{2}}(\mathbb{R}^{3}_{x})$ is not included in $L^\infty$. Thus, it seems more natural to
choose  homogeneous Besov space $\dot{B}^{\f{3}{2}}_{2,1}$ for
the density $c$ and the magnetic field $H$ since $\dot{B}^{\f{3}{2}}_{2,1}$ is continuously
embedded into $L^\infty$. As in \cite{CMZ2,Dan2,Dan1,Dan3}, the different dissipative
mechanisms of low frequencies and high frequencies inspire us to
deal with $c$ and $H$ in $\tilde{B}^{\f12,\f32}_{2,1}=\dot{B}^{\f12}_{2,1}\cap\dot{B}^{\f32}_{2,1}$. However,  we  can not obtain the  desired bounds directly  in critical regularity framework.
In particular, there is a difficulty coming from the convection term $u\cdot\na c$ in transport equations  in high frequencies, as  one derivative loss about the function $c$ will appear no matter how smooth is $u$ if it is viewed as a perturbation term.  To overcome
the difficulty, we need establish a uniform estimate for a mixed hyperbolic-parabolic linear system
with  convection terms (See Proposition \ref{prop3.1}). Combining  with some nonlinear  estimates and  the standard continuity arguments, we  obtain global existence and uniqueness of strong solutions to the system \eqref{equ:CMHD} in critical Besov spaces if the initial density is close to a positive constant. Next,  one may wonder  how global strong solutions constructed above look like for large time. Under a suitable additional condition involving only the low frequencies of the  data and in $L^{2}$-critical regularity framework,  we  exhibit the optimal time decay rates  for  the  constructed global strong
solutions. In this part, our main ideas are based on an application of Fourier analysis to a linearized parabolic-hyperbolic system and on a refined time-weighted energy functional.  In  low frequencies, making good use of the decay estimates of
Green's function for the linearized system  and  Duhamel's  principle, one can obtain the
desired estimates.  In  high frequencies,  we can deal with  the estimates of the nonlinear terms in the
system employing  the Fourier localization method, the energy method
and commutator estimates. Finally, in order to close the energy estimates, we exploit some  decay estimates with gain of
regularity  for the high frequencies of $\nabla u, \nabla^{2}H$.

Now we state our main results as follows:
\begin{Theorem}\label{th:main1}  Assume that $(\rho _{0}-1,u_{0},H_{0})\in \tilde{B}^{\f{1}{2},\f{3}{2}}_{2,1}\times
\dot{B}^{\f{1}{2}}_{2,1}\times \tilde{B}^{\f{1}{2},\f{3}{2}}_{2,1}$ and  with
no loss of generality that $\bar{\rho}=1$ and $p'(1)=1$. Then there exists a constant $\eta>0$ such that if
\begin{equation}
\|\rho _{0}-1\|_{\tilde{B}^{\f{1}{2},\f{3}{2}}_{2,1}}+\|u_{0}\|_{
\dot{B}^{\f{1}{2}}_{2,1}}+\|H_{0}\|_{\tilde{B}^{\f{1}{2},\f{3}{2}}_{2,1}}\leq \eta,  \label{1Hn/2}
\end{equation}%
then the Cauchy  problem \eqref{equ:CMHD}-\eqref{initial value} admits a unique global solution $(\rho
-1,u,H)$ satisfying that for all $t\geq 0$,
\begin{equation}
\begin{split}
\label{1.1}
\,X(t)\eqdefa&\|(\rho-1,u,H)\|_{\wt L^\infty_t(\dot B^{\frac{1}{2}}_{2,1})}^\ell
+\|(\rho-1,u,H)\|_{L^1_t(\dot B^{\frac{5}{2}}_{2,1})}^\ell+\|u\|_{\wt L^\infty_t(\dot B^{\frac{1}{2}}_{2,1})}^h
+\|(\rho-1,H)\|_{\wt L^\infty_t(\dot B^{\frac{3}{2}}_{2,1})}^h
\\
&\quad+\|u\|_{L^1_t(\dot B^{\frac{5}{2}}_{2,1})}^h
+\|\rho-1\|_{L^1_t(\dot B^{\frac{3}{2}}_{2,1})}^h
+\|H\|_{L^1_t(\dot B^{\frac{7}{2}}_{2,1})}^h
\\&\quad \lesssim \|\rho _{0}-1\|_{\tilde{B}^{\f{1}{2},\f{3}{2}}_{2,1}}+\|u_{0}\|_{
\dot{B}^{\f{1}{2}}_{2,1}}+\|H_{0}\|_{\tilde{B}^{\f{1}{2},\f{3}{2}}_{2,1}}.
\end{split}
\end{equation}
\end{Theorem}
\begin{Theorem}\label{th:decay}
Let the data $(\rho _{0}-1,u_0,H_0)$ satisfy the assumptions of Theorem \ref{th:main1}. Denote $\langle \tau\rangle\eqdefa\sqrt{1+\tau^2}$
and $\alpha\eqdefa\frac{5}{2}-\varepsilon$ with $\varepsilon>0$ arbitrarily small.
 There exists a positive constant $c$ so that if in addition
\begin{equation}\label{eq:D0}
D_0\eqdefa\|(\rho_0-1,u_0,H_0)\|^{\ell}_{\dot B^{-\frac 32}_{2,\infty}}\leq c,
\end{equation}
then  the global solution $(\rho-1,u,H)$
given by Theorem \ref{th:main1} satisfies for all $t\geq0,$
\begin{equation}
\label{1.8}
D(t)\leq C\bigl(D_0+\|(\nabla \rho_{0}, u_{0},\nabla H_{0})\|^h_{\dot B^{\frac 12}_{2,1}}\bigr)
\end{equation}
with
\begin{align*}
D(t)\eqdefa
&\sup_{s\in(-\frac 32,2]}\|\langle\tau\rangle^{\frac 34+\frac s2}(\rho-1,u,H)\|_{L^\infty_t(\dot B^s_{2,1})}^\ell
+\|\langle\tau\rangle^{\alpha}(\nabla \rho,u,\nabla H)\|_{\wt L^\infty_t(\dot B^{\frac 12}_{2,1})}^h\\
&+\|\tau\nabla  u\|_{\wt L^\infty_t(\dot B^{\frac{3}{2}}_{2,1})}^h
+\|\tau\nabla^{2}  H\|_{\wt L^\infty_t(\dot B^{\frac{3}{2}}_{2,1})}^h.
\end{align*}
\end{Theorem}
\begin{Remark} In Theorem \ref{th:main1}, we extend the existence result of 3D incompressible Hall-MHD system in critical Besov spaces in \cite{chae3} to 3D compressible case.
\end{Remark}
\begin{Remark}\label{1.2} In Theorem \ref{th:decay}, we obtain the optimal decay rates for the 3D compressible Hall-MHD equations
\eqref{equ:CMHD}  in critical regularity framework. Additionally,
the regularity  index $s$ can take both negative and nonnegative values, rather than only nonnegative integers, which   improves the classical decay results  in high Sobolev regularity, such as \cite{fan3,GY,Xu3}. In fact, for the solution $(\rho-1,u,H)$  constructed in  Theorem \ref{th:decay}, employing homogeneous Littlewood-Paley decomposition
 for $\rho-1$, we have
 $$\Lambda^s (\rho-1)=\sum_{q\in \mathbb{Z}}\dot{ \Delta}_q\Lambda^s (\rho-1).$$  Thus
$$\|\Lambda^s (\rho-1)\|_{L^2}\lesssim  \sum_{q\in \mathbb{Z}} \|\Delta_q\Lambda^s (\rho-1)\|_{L^2}=\|\Lambda^s (\rho-1)\|_{\dot B^0_{2,1}}.$$
Based on the Bernstein inequalities and the low-high frequencies decomposition, we may write
$$
\sup_{t\in[0,T]} \langle t\rangle^{\frac {3/2+s}2}\|\Lambda^s (\rho-1)\|_{\dot B^0_{2,1}}\lesssim
 \|\langle t\rangle^{\frac {3/2+s}2}(\rho-1)\|_{L^\infty_T(\dot B^s_{2,1})}^\ell
 +  \|\langle t\rangle^{\frac {3/2+s}2}(\rho-1)\|_{L^\infty_T(\dot B^s_{2,1})}^h.
$$
If follows from Inequality \eqref{1.8} and  definitions of $D(t)$ and $\alpha$ that
$$
  \|\langle t\rangle^{\frac {3/2+s}2}\rho-1\|_{L^\infty_T(\dot B^s_{2,1})}^\ell\lesssim D_0+\|(\nabla \rho_{0}, u_{0},\nabla H_{0})\|^h_{\dot B^{\frac 12}_{2,1}}
\quad\hbox{if }\ \  -3/2<s\leq 2
$$
and that, because we have $\alpha\geq\frac{3/2+s}2$ for all $s\leq3/2,$
$$
  \|\langle t\rangle^{\frac {3/2+s}2}\rho-1\|_{L^\infty_T(\dot B^s_{2,1})}^h
  \lesssim D_0+\|(\nabla \rho_{0}, u_{0},\nabla H_{0})\|^h_{\dot B^{\frac 12}_{2,1}}\quad\hbox{if }\ \  s\leq 3/2.$$
This yields the following desired result for $\rho-1$
$$\|\Lambda^{s}(\rho-1)\|_{L^2}
\leq C\bigl(D_0+\|(\nabla \rho_{0}, u_{0},\nabla H_{0})\|^h_{\dot B^{\frac 12}_{2,1}}\bigr)\langle t\rangle^{-\frac {3/2+s}2}
\quad  \hbox{ if } \ \ -3/2<s\leq3/2,$$  where the fractional derivative
 operator $\Lambda^{\ell}$ is defined by $\Lambda^{\ell}f\triangleq\mathcal{F}^{-1}(|\cdot|^{\ell}\mathcal{F}f)$. Similarly, we have
$$\displaylines{
\|\Lambda^{s}u\|_{L^2}
\leq C\bigl(D_0+\|(\nabla \rho_{0}, u_{0},\nabla H_{0})\|^h_{\dot B^{\frac 12}_{2,1}}\bigr)\langle t\rangle^{-\frac {3/2+s}2}
\quad \hbox{ if } \ \ -3/2<s\leq1/2,\cr
\|\Lambda^{s}H\|_{L^2}
\leq C\bigl(D_0+\|(\nabla \rho_{0}, u_{0},\nabla H_{0})\|^h_{\dot B^{\frac 12}_{2,1}}\bigr)\langle t\rangle^{-\frac {3/2+s}2}
\quad  \hbox{ if } \ \ -3/2<s\leq3/2.}
$$
In particular, taking $s=0$ leads back to  the standard optimal   $L^{1}$-$L^{2}$ decay rate of $(\rho-1,u,H)$ as in \cite{fan3,GY}.
 \end{Remark}
\begin{Remark}\label{1.2}Due to the embedding $L^{1}(\mathbb{R}^3)\hookrightarrow \dot{B}^{-\frac 32}_{2,\infty}(\mathbb{R}^3)$,  our results in Theorem \ref{th:decay} extend the known conclusions in \cite{fan3,GY}. In particular, our condition involves only the low frequencies of the data and is based on $L^{2}(\mathbb{R}^{3}_{x})$-norm framework.
\end{Remark}
\noindent{\bf Notaions.}  We assume $C$ be
a positive generic constant throughout this paper that may vary at
different places and
denote  $A\le CB$ by  $A\lesssim B$.
\par
\section{ Littlewood-Paley theory and some useful lemmas }
\ \ \ Let us introduce the Littlewood-Paley decomposition.
Choose a radial function  $\varphi \in {\cS}(\mathbb{R}^N)$
supported in ${\cC}=\{\xi\in\mathbb{R}^N,\,
\frac{3}{4}\le|\xi|\le\frac{8}{3}\}$ such that \beno \sum_{q\in
\mathbb{Z}}\varphi(2^{-q}\xi)=1 \quad \textrm{for all}\,\,\xi\neq 0.
\eeno The homogeneous frequency localization operators
$\dot{\Delta}_q$ and $\dot{S}_q$ are defined by
\begin{align}
\dot{\Delta}_qf=\varphi(2^{-q}D)f,\quad \dot{S}_qf=\sum_{k\le
q-1}\dot{\Delta}_kf\quad\mbox{for}\quad q\in \mathbb{Z}. \nonumber
\end{align}
With our choice of $\varphi$, one can easily verify that
\begin{equation*}\begin{split}
&\dot{\Delta}_q\dot{\Delta}_kf=0\quad \textrm{if}\quad|q-k|\ge
2\quad \textrm{and}
\quad \\
&\dot{\Delta}_q(\dot{S}_{k-1}f\dot{\Delta}_k f)=0\quad
\textrm{if}\quad|q-k|\ge 5.
\end{split}\end{equation*}
 We denote the space ${\cZ'}(\mathbb{R}^N)$ by the dual space of
${\cZ}(\mathbb{R}^N)=\{f\in {\cS}(\mathbb{R}^N);\,D^\alpha
\hat{f}(0)=0; \forall\alpha\in\mathbb{ N}^N \,\mbox{multi-index}\}$,
it also can be identified by the quotient space of
${\cS'}(\mathbb{R}^N)/{\cP}$ with the polynomials space ${\cP}$. The
formal equality \beno f=\sum_{q\in\mathbb{Z}}\dot{\Delta}_qf \eeno
holds true for $f\in {\cZ'}(\mathbb{R}^N)$ and is called the
homogeneous Littlewood-Paley decomposition.

The following Bernstein's inequalities will be frequently used.
\begin{Lemma}\cite{Che-book}\label{Lem:Bernstein}
Let $1\le p_{1}\le p_{2}\le+\infty$. Assume that $f\in L^{p_{1}}(\mathbb{R}^N)$,
then for any $\gamma\in(\mathbb{N}\cup\{0\})^N$, there exist
constants $C_1$, $C_2$ independent of $f$, $q$ such that \beno
&&{\rm supp}\hat f\subseteq \{|\xi|\le A_02^{q}\}\Rightarrow
\|\partial^\gamma f\|_{p_2}\le C_12^{q{|\gamma|}+q
N(\frac{1}{p_1}-\frac{1}{p_2})}\|f\|_{p_1},
\\
&&{\rm supp}\hat f\subseteq \{A_12^{q}\le|\xi|\le
A_22^{q}\}\Rightarrow \|f\|_{p_1}\le
C_22^{-q|\gamma|}\sup_{|\beta|=|\gamma|}\|\partial^\beta f\|_{p_1}.
\eeno
\end{Lemma}
\begin{Definition} Let $s\in \mathbb{R}$, $1\le p,
r\le+\infty$. The homogeneous Besov space $\dot{B}^{s}_{p,r}$ is
defined by
$$\dot{B}^{s}_{p,r}=\Big\{f\in {\cZ'}(\mathbb{R}^N):\,\|f\|_{\dot{B}^{s}_{p,r}}<+\infty\Big\},$$
where \beno \|f\|_{\dot{B}^{s}_{p,r}}\eqdefa \Big\|2^{qs}
\|\dot{\Delta}_qf(t)\|_{p}\Big\|_{\ell^r}.\eeno
\end{Definition}
\begin{Remark}\label{2.3}
Some properties about the  Besov spaces are as follows
\begin{itemize}
\item\,\, Derivation: $$\|f\|_{\dot{B}^{s}_{2,1}}\approx\|\nabla f\|_{\dot{B}^{s-1}_{2,1}};$$
\item\,\, Algebraic properties: for $s > 0$, $\dot{B}^{s}_{2,1} \cap L^{\infty}
$ is an algebra;
\item\,\,Interpolation: for
$s_1, s_2\in\mathbb{R}$ and $\theta\in[0,1]$,
we have $$\|f\|_{\dot{B}^{\theta s_1+(1-\theta)s_2}_{2,1}}\le
\|f\|^\theta_{\dot{B}^{s_1}_{2,1}}\|f\|^{(1-\theta)}_{\dot{B}^{s_2}_{2,1}}.$$
\end{itemize}
\end{Remark}
\begin{Definition} Let $s\in \mathbb{R}$, $1\le p,\rho,
r\le+\infty$. The homogeneous space-time  Besov space
$L^\rho_{T}(\dot{B}_{p,  r}^s)$ is defined by
$$L^\rho_{T}(\dot{B}_{p,  r}^s)=\Big\{f\in \mathbb{R}_{+}\times{\cZ'}(\mathbb{R}^N):\,\|f\|_{L^\rho_{T}(\dot{B}_{p,  r}^s)} <+\infty\Big\},$$ where $$
\|f\|_{L^\rho_{T}(\dot{B}_{p,  r}^s)}\eqdefa \Big\|\big\| 2^{qs}
\|\dot{\Delta}_q f\|_{L^p}\big\|_{\ell ^{r}}\Big\|_{L^\rho_{T}}.$$
\end{Definition}
We next introduce the Besov-Chemin-Lerner space
$\widetilde{L}^q_T(\dot{B}^{s}_{p,r})$ which is initiated in
\cite{Che-Ler}.
\begin{Definition}Let $s\in \mathbb{R}$, $1\le
p,q,r\le+\infty$, $0<T\le+\infty$. The space
$\widetilde{L}^q_T(\dot{B}^s_{p,r})$ is defined by
$$\widetilde{L}^q_T(\dot{B}^s_{p,r})=\Big\{f\in \mathbb{R}_{+}\times{\cZ'}(\mathbb{R}^N):\,\|f\|_{\widetilde{L}^q_T(\dot{B}^{s}_{p,r})}<+\infty\Big\},$$
where
$$\|f\|_{\widetilde{L}^q_T(\dot{B}^{s}_{p,r})}\eqdefa \Bigl\|2^{qs}
\|\dot{\Delta}_qf(t)\|_{L^q(0,T;L^p)}\Bigr\|_{\ell^r}.$$
\end{Definition}
Obviously, $
\widetilde{L}^1_T(\dot{B}^s_{p,1})=L^1_T(\dot{B}^s_{p,1}). $ By  a
direct application  of  Minkowski's inequality, we have the
following relations between these spaces
\begin{equation*}
L^\rho_{T}(\dot{B}_{p,r}^s)\hookrightarrow\widetilde
L^\rho_{T}(\dot{B}_{p,r}^s),\,\textnormal{if}\quad  r\geq
\rho,\end{equation*}
\begin{equation*}
\widetilde L^\rho_{T}(\dot{B}_{p,r}^s)\hookrightarrow
L^\rho_{T}(\dot{B}_{p,r}^s),\, \textnormal{if}\quad \rho\geq r.
\end{equation*}
To deal with functions with different regularities for high
frequencies and low frequencies, motivated by \cite{Dan2,Dan1}, it
is more effective to work in \textit{hybrid Besov spaces}. We remark
that using hybrid Besov spaces has been crucial for proving global
well-posedness for compressible systems in critical spaces (see
\cite{CMZ2,Dan2,Dan1}).
\begin{Definition}
Let $s,t\in \mathbb{R}$. We set
$$\|f\|_{\tilde{B}^{s,t}_{2,1}}=\sum_{q\le
0}2^{qs}\|\dot{\D}_qf\|_{L^2}+\sum_{q>0}2^{qt}\|\dot{\D}_qf\|_{L^2}.$$
For $m=-\left[\f{N}{2}+1-s\right]$, we define
\begin{eqnarray}
\tilde{B}^{s,t}_{2,1}&=&\left\{f\in\mathcal{S}'(\mathbb{R}^N):
\|f\|_{\tilde{B}^{s,t}_{2,1}}<\infty\right\},\qquad \textrm{
if}\ m<0,\nonumber\\
\tilde{B}^{s,t}_{2,1}&=&\left\{f\in\mathcal{S}'(\mathbb{R}^N)/\mathcal{P}_m:
\|f\|_{\tilde{B}^{s,t}_{2,1}}<\infty\right\},\ \textrm{ if }m\ge
0.\nonumber
\end{eqnarray}
\end{Definition}
\begin{Remark}\label{2.5}
Some properties about the hybrid Besov spaces are as follows
\begin{itemize}
\item\,\, $\tilde{B}^{s,s}_{2,1}=\dot{B}^s_{2,1}$;
\item\,\, If $s\le t$, then $\tilde{B}^{s,t}_{2,1}=\dot{B}^s_{2,1}\cap \dot{B}^t_{2,1}$. Otherwise,
$\tilde{B}^{s,t}_{2,1}=\dot{B}^s_{2,1}+\dot{B}^t_{2,1}$. In
particular, $\tilde{B}^{s,\f{N}{2}}_{2,1}\hookrightarrow L^\infty$
as $s\le \f{N}{2}$;
\item\,\,Interpolation: for
$s_1, s_2, t_1, t_2\in\mathbb{R}$ and $\theta\in[0,1]$,
we have $$\|f\|_{\tilde{B}^{\theta s_1+(1-\theta)s_2,\, \theta t_1+(1-\theta)t_2}_{2,1}}\le
\|f\|^\theta_{\tilde{B}^{s_1,t_1}_{2,1}}\|f\|^{(1-\theta)}_{\tilde{B}^{s_2,t_2}_{2,1}};$$
\item\,\, If $s_1\le s_2$ and $t_1\ge t_2$, then
$\tilde{B}^{s_1,t_1}_{2,1}\hookrightarrow\tilde{B}^{s_2,t_2}_{2,1}$.
\end{itemize}
\end{Remark}
We have the following properties for the product in  Besov spaces and  hybrid Besov spaces.
\begin{Proposition}\cite{Dan5}\label{p26}
 For all $1\leq r,p, p_1, p_2\leq+\infty$,  there exists a positive universal
 constant such that
$$\|fg\|_{\dot{B}^{s}_{p,r}}\lesssim
\|f\|_{L^\infty}\|g\|_{\dot{B}^{s}_{p,r}}+\|g\|_{L^\infty}\|f\|_{\dot{B}^{s}_{p,r}},
\quad \text{if}\quad s>0;$$
$$\|fg\|_{\dot{B}^{s_1+s_2-\frac{N}{p}}_{p,r}}\lesssim
\|f\|_{\dot{B}^{s_1}_{p,r}}\|g\|_{\dot{B}^{s_2}_{p,\infty}}, \quad
\text{if}\quad s_1,s_2<\frac{N}{p},\quad \text{and}\quad
s_1+s_2>0;$$
$$\|fg\|_{\dot{B}^{s}_{p,r}}\lesssim
\|f\|_{\dot{B}^{s}_{p,r}}\|g\|_{\dot{B}^{\frac{N}{p}}_{p,\infty}\cap
L^{\infty}}, \quad \text{if}\quad |s|<\frac{N}{p};$$
$$\|fg\|_{\dot{B}^s_{2,1}}\lesssim \|f\|_{\dot{B}^{d/2}_{2,1}}\|g\|_{\dot{B}^s_{2,1}}, \quad \text{if}\quad
s\in (-d/2,d/2].$$
\end{Proposition}
\begin{Proposition}\cite{Dan3}\label{p25}
 For all $s_1, s_2>0$,  there exists a positive universal
 constant such that
$$\|fg\|_{\tilde{B}^{s_1,s_2}_{2,1}}\lesssim
\|f\|_{L^\infty}\|g\|_{\tilde{B}^{s_1,s_2}_{2,1}}+\|g\|_{L^\infty}\|f\|_{\tilde{B}^{s_1,s_2}_{2,1}}.$$
For all $s_1, s_2\le\f{N}{2}$ such that $\min\{s_1+t_1,
s_2+t_2\}>0$,  there exists a positive universal
 constant such that
$$\|fg\|_{\tilde{B}^{s_1+t_1-\f{N}{2},s_2+t_2-\f{N}{2}}_{2,1}}\lesssim
\|f\|_{\tilde{B}^{s_1,s_2}_{2,1}}\|g\|_{\tilde{B}^{t_1,t_2}_{2,1}}.$$
\end{Proposition}
For the composition of functions, we have the following estimates.
\begin{Proposition}\cite{Dan1}\label{p27}
Let $s>0$ and $u\in \dot{B}^{s}_{2.1}\cap L^\infty$.

\mbox{(i) } If $F\in W_{loc}^{[s]+2,\infty}(\mathbb{R}^N)$ with
$F(0)=0$, then $F(u)\in \dot{B}^{s}_{2,1}$. Moreover, there exists a
function of one variable $C_0$ depending only on $s$ and $F$, and
such that
\begin{align*}
    \|F(u)\|_{\dot{B}^{s}_{2,1}}\leq
    C_0(\|u\|_{L^\infty})\|u\|_{\dot{B}^{s}_{2,1}}.
\end{align*}

\mbox{ (ii)} If $u,\, v\in \dot{B}^{\frac{N}{2}}_{2,1}$, $(v-u)\in
\dot{B}^{s}_{2,1}$ for $s\in(-\frac{N}{2},\frac{N}{2}]$ and $G\in
W_{loc}^{[\frac{N}{2}]+3,\infty}(\mathbb{R}^N)$ satisfies $G'(0)=0$,
then $G(v)-G(u)\in \dot{B}^{s}_{2,1}$ and there exists a function of
two variables $C$ depending only on $s$, $N$ and $G$, and such that
\begin{align*}
  \|G(v)-G(u)\|_{\dot{B}^{s}_{2,1}}\leq C\big(\|u\|_{L^\infty},
    \|v\|_{L^\infty}\big)\big(\|u\|_{\dot{B}^{\frac{N}{2}}_{2,1}}+\|v\|_{\dot{B}^{\frac{N}{2}}_{2,1}}\big)
   \|v-u\|_{\dot{B}^{s}_{2,1}}.
\end{align*}
\end{Proposition}
Throughout this paper, the following estimates for the convection
terms arising in the linearized systems will be used frequently.
\begin{Proposition}\cite{Dan1}\label{p29}
    Let $F$ be an homogeneous smooth function of degree $m$. Suppose that $-N/2<s_1,t_1,s_2,t_2\leq 1+N/2$. The following two estimates hold
    \begin{align*}
    &|(F(D)\dot{\Delta}_{q}(v\cdot \nabla a),F(D)\dot{\Delta}_{q} a)|\\
    &\qquad\qquad\leq C \gamma_q2^{-q(\phi^{s_1,s_2}(q)-m)}\|v\|_{\dot{B}^{\frac{N}{2}+1}_{2,1}} \|a\|_{\tilde{B}^{s_,s_2}_{2,1}}\|F(D)\dot{\Delta}_{q} a\|_{2},\\
    &|(F(D)\dot{\Delta}_{q}(v\cdot\nabla a),\dot{\Delta}_{q} b)+(\dot{\Delta}_{q}(v\cdot\nabla b),F(D)\dot{\Delta}_{q} a)|\\
    &\qquad\qquad\leq C \gamma_q\|v\|_{\dot{B}^{\frac{N}{2}+1}_{2,1}}\times\big( 2^{-q\phi^{t_1,t_2}(q)}\|F(D)\dot{\Delta}_{q} a\|_{2}\|b\|_{\tilde{B}^{t_,t_2}_{2,1}}\\
    &\qquad\qquad\qquad+2^{-q(\phi^{s_1,s_2}(q)-m)}\|a\|_{\tilde{B}^{s_,s_2}_{2,1}}
    \|\dot{\Delta}_{q} b\|_{2}\big),
    \end{align*}
where $(\cdot,\cdot)$ denotes the $2$-inner product,  $\sum_{q\in
\mathbb{Z}} \gamma_q\leq 1$ and the operator $F(D)$ is defined by
$F(D)f:=\mathcal{F}^{-1} F(\xi)\mathcal{F} f$,
$\phi^{\alpha,\beta}(q)$ is the following characteristic function on
$\mathbb{Z}$
\begin{align*}
     \phi^{\alpha,\beta}(q)=\left\{
     \begin{array}{ll}
        \alpha,\quad &\text{if }\quad  q\leq 0,\\
        \beta, \quad & \text{if }\quad q\geq 1.
     \end{array}
     \right.
\end{align*}
\end{Proposition}
\begin{Proposition}\cite{Dan2}
\label{Pro:1}
Let $1\leq p,p_{1}\leq\infty$, $1\leq r\leq\infty$ and $\sigma\in\mathbb{R}$.
There exists a constant $C>0$ depending only on $\sigma$ such that for all $q\in \mathbb{Z}$, we have
$$\|[v\cdot\nabla, \partial_{\ell}\dot{\Delta}_{q}]a\|_{L^{p}}
\leq Cc_{q}2^{-q(\sigma-1)}
\|\nabla v\|_{\dot B^{\frac{d}{p_{1}}}_{p_{1},1}}
\|\nabla a\|_{\dot B^{\sigma-1}_{p,1}}, \quad \text{for}\quad -\min(\frac{d}{p_{1}},\frac{d}{p^{\prime}})<\sigma
\leq1+\min(\frac{d}{p_{1}},\frac{d}{p}),$$
$$\|[v\cdot\nabla,\dot{\Delta}_{q}]a\|_{L^{p}}
\leq Cc_{q}2^{-q\sigma}
\|\nabla v\|_{\dot B^{\frac{d}{p_{1}}}_{p_{1},\infty}\cap L^{\infty}}
\|a\|_{\dot B^{\sigma}_{p,1}}, \quad \text{for}\quad -\min(\frac{d}{p_{1}},\frac{d}{p^{\prime}})<\sigma
<1+\frac{d}{p_{1}},$$
where the commutator $[\cdot,\cdot]$ is defined by $[f,g]=fg-gf$ and $(c_{j})_{j\in \mathbb{Z}}$ denotes a
sequence such that $\sum_{q\in
\mathbb{Z}} c_{q}\leq 1$.
\end{Proposition}
\begin{Proposition}\cite{RD}
\label{Pro:3}
Let $\sigma\in\mathbb{R}, (p,r)\in[1,\infty]^{2}$ and $1\leq \rho_{2}\leq \rho_{1}\leq\infty$. Let $u$ satisfy
\begin{align}\label{eq:heat}
\left\{
\begin{aligned}
&\partial_{t}u-\mu\Delta u=f,\\
&u\mid_{t=0}=u_{0}. \end{aligned} \right.
\end{align}
Then for all $T>0$ the following a priori estimates is fulfilled
\begin{equation}\label{eq:heat1}\mu^{\frac{1}{\rho_{1}}}\|u\|_{\wt L^{\rho_{1}}_T(\dot B^{\sigma+\frac 2\rho_{1}}_{p,r})}\lesssim\|u_{0}\|_{\dot B^{\sigma}_{p,r}}
+\mu^{\frac{1}{\rho_{2}}-1}\|f\|_{\wt L^{\rho_{2}}_T(\dot B^{\sigma-2+\frac 2\rho_{2}}_{p,r})}.
\end{equation}
\end{Proposition}
\begin{Remark}\label{2.14}
The solutions to the following Lam\'{e} system
\begin{align*}\label{eq:heat}
\left\{
\begin{aligned}
&\partial_tu-\cA u=f,\\
&u|_{t=0}=u_0,
\end{aligned} \right.
\end{align*}
are also fulfill \eqref{eq:heat1}.
\end{Remark}
We finish this subsection by listing an elementary but useful
inequality.
\begin{Lemma}\cite{MN}\label{lemma2.13}\ (a)\  Let $r_1,r_2>0$ satisfy $\max\{r_1,r_2\}>1$. Then
$$\int_0^t(1+t-\tau)^{-r_1}(1+\tau)^{-r_2}d\tau\leq C(r_1,r_2)(1+t)^{-\min\{r_1,r_2\}}.$$
\ (b)\  Let $r_1,r_2>0$ and $f\in L^{1}(0,+\infty)$. Then
$$\int_0^t(1+t-\tau)^{-r_1}(1+\tau)^{-r_2}f(\tau)d\tau\leq C(r_1,r_2)(1+t)^{-\min\{r_1,r_2\}}\int_{0}^{t}|f|d\tau.$$
\end{Lemma}
\par
\section{Reformulation of the Original System \eqref{equ:CMHD} and {A priori} estimates for linearized  system with convection terms}

\subsection{Reformulation of the Original System \eqref{equ:CMHD}}
\ \  We first reformulate the original system
\eqref{equ:CMHD} into a different form. For the magnetic field $H$,
we have the following identities:
$$\nabla(|H|^{2})=2(H\cdot\nabla)H+2(\nabla\times H)\times H,$$
$$\nabla\times(\nabla\times H)=\nabla\Dv H-\Delta H,$$
and  \begin{equation*}
\begin{split} \nabla\times(u\times H)&=u(\Dv H)-H(\Dv u)+H\cdot\nabla
u-u\cdot\nabla H\\&=-H(\Dv u)+H\cdot\nabla
u-u\cdot\nabla H\end{split}
\end{equation*} with $\Dv H=0$.

We will denote that $c=\rho-1$. Then, in term of the new variables $(c, u,
H)$,  system (\ref{equ:CMHD})-(\ref{initial value}) becomes
\begin{equation}\label{equ:heatflow-new}
\left\{
\begin{aligned}{}
&\p_tc+u\cdot\na c+\dv u=f,\\
&\p_t u+u\cdot\na u-\mathcal{A}u+\nabla c=g,\\
&\p_tH+u\cdot\nabla H-\Delta H=h, \\
&\textrm{div}H=0,\\
&(c,\,u,\,H)|_{t=0}=\big(c_0,\,u_0,\,H_0\big),
\end{aligned}
\right.
\end{equation}
where
\begin{equation*}
\begin{split}f&=-c\dv u,\\
g&=-L_1(c)\mathcal{A}u+ L_2(c)\na
c-L_3(c)\Big(\frac{1}{2}\nabla |H|^{2}-H\cdot\nabla H\Big),\\
h&=-H(\dv u)+H\cdot\nabla u -\nabla\times \Big(L_3(c)(\nabla\times H)\times H\Big),
\end{split}
\end{equation*} with
$$ \mathcal{A}u=\mu\Delta u+(\lambda+\mu)\na
\dv u, \quad L_1(c)=\frac{c}{1+c}, \quad
L_2(c)=\frac{p'(1+c)}{1+c}-1, \quad L_3(c)=\frac{1}{c+1}.$$
\subsection{{A priori} estimates for linearized  system with convection terms}
\ \  Next, we investigate some \emph{a priori} estimates for the following linearized  system with convection terms
\begin{equation}\label{equ:heatflow-new-11}
\left\{
\begin{aligned}{}
&\p_tc+v\cdot\na c+\dv u=f,\\
&\p_t u+v\cdot\na u-\mu\Delta u+(\lambda+\mu)\na
\dv u+\nabla c=g,\\
&\p_tH+v\cdot\nabla H- \Delta H=h, \\
&(c,\,u,\,H)|_{t=0}=\big(c_0,\,u_0,\,H_0\big).
\end{aligned}
\right.
\end{equation}
We prove the following
proposition and  show a uniform estimate for a mixed hyperbolic-parabolic linear system
with  convection terms. What is crucial in this work is the smoothing effect on the velocity $u$ and magnetic fields $H$ and a $L^{1}$ decay on
$\r-\bar{\r}$ (this plays a key role to control the pressure term).
\begin{Proposition}\label{prop3.1}
 Denote
\begin{equation*}
V(t):=\int_0^t\|v(\tau)\|_{\dot{B}^{\f{5}{2}}_{2,1}}d\tau.
\label{VV}
\end{equation*}
Let $(c,\,u,\,H)$ be a solution of the system
\eqref{equ:heatflow-new-11} on $[0,t)$. Then the following estimates
hold  for  $t\in[0,T)$
\begin{equation}\label{4.21}
\begin{split}
&\bigl\|c\bigr\|_{\widetilde{L}^{\infty}([0,t];\tilde{B}^{\f{1}{2},\f{3}{2}}_{2,1})}+\bigl\|u\bigr\|_{\widetilde{L}^{\infty}([0,t];
\dot{B}^{\f{1}{2}}_{2,1})}+\bigl\|H\bigr\|_{\widetilde{L}^{\infty}([0,t];\tilde{B}^{\f{1}{2},\f{3}{2}}_{2,1})}
\\&\qquad+\int_0^t\big\|c(\tau)\big\|_{\tilde{B}^{\f{5}{2},\f{3}{2}}_{2,1}}d\tau+\int_0^t\big\|u(\tau)\big\|_{
\dot{B}^{\f{5}{2}}_{2,1}}d\tau+\int_0^t\big\|H(\tau)\big\|_{\tilde{B}^{\f{5}{2},\f{7}{2}}_{2,1}}d\tau
\\&\quad \lesssim e^{CV(t)}\Big( \|c_{0}\|_{\tilde{B}^{\f{1}{2},\f{3}{2}}_{2,1}}+\|u_{0}\|_{
\dot{B}^{\f{1}{2}}_{2,1}}+\|H_{0}\|_{\tilde{B}^{\f{1}{2},\f{3}{2}}_{2,1}}+\int_0^t\big\|f(\tau)\big\|_{\tilde{B}^{\f{1}{2},\f{3}{2}}_{2,1}}d\tau
\\&\qquad+\int_0^t\big\|g(\tau)\big\|_{\dot{B}^{\f{1}{2}}_{2,1}}d\tau+\int_0^t\big\|h(\tau)\big\|_{\tilde{B}^{\f{1}{2},\f{3}{2}}_{2,1}}d\tau\Big).
\end{split}
\end{equation}
\end{Proposition}
\noindent{\bf Proof.} To prove the proposition, we first localize
the system \eqref{equ:heatflow-new-11}  according to the homogeneous Littlewood-Paley decomposition. Then
each dyadic block can be estimated by using energy method. Applying the operator $\dot{\Delta}_q$ to the
system \eqref{equ:heatflow-new-11}, we deduce that
$(\dot{\Delta}_qc,\dot{\Delta}_qu,\dot{\Delta}_qH
)$ satisfies
\begin{equation}\label{equ:heatflow-new31}
\left\{
\begin{aligned}{}
&\p_t\dot{ \Delta}_qc+\dot{ \Delta}_q(v\cdot\na c)+\dv \dot{ \Delta}_qu=\dot{ \Delta}_qf,\\
&\p_t \dot{ \Delta}_qu+\dot{ \Delta}_q(v\cdot\nabla u)-\mu\Delta \dot{ \Delta}_qu+(\lambda+\mu)\na
\dv \dot{ \Delta}_qu+\nabla\dot{ \Delta}_qc=\dot{ \Delta}_qg,\\
&\p_t \dot{ \Delta}_qH+\dot{ \Delta}_q(v\cdot\nabla H)-\Delta \dot{ \Delta}_qH=\dot{ \Delta}_qh.
\end{aligned}
\right.
\end{equation}
Taking the $L^2$-scalar product of the
first equation of \eqref{equ:heatflow-new3} with $\dot{ \Delta}_qc$ and  $(-\Delta\dot{\Delta}_q c)$,
the second equation with $\dot{\Delta}_qu$,  the third equation with
$\dot{\Delta}_qH$ and $\Delta\dot{\Delta}_qH$ respectively,  we obtain the
following five identities
\begin{equation}
\label{11}\begin{split}
\f{1}{2}\f{d}{dt}\|\dot{\Delta}_qc\|_{L^2}^2+(\dot{\Delta}_q(v\cdot\nabla c)|\dot{\Delta}_qc)+(\dv\dot{\Delta}_qu|\dot{\Delta}_qc)
=(\dot{\Delta}_qf|\dot{\Delta}_qc),
\end{split}
\end{equation}
\begin{equation}
\label{12}\begin{split}
\f{1}{2}\f{d}{dt}\|\nabla\dot{\Delta}_qc\|_{L^2}^2+(\dot{\Delta}_q(v\cdot\nabla c)|-\Delta\dot{\Delta}_qc)+(\dv\dot{\Delta}_qu|-\Delta\dot{\Delta}_qc)
=(\dot{\Delta}_qf|-\Delta\dot{\Delta}_qc),
\end{split}
\end{equation}
\begin{equation}
\label{13}\begin{split}
\f{1}{2}\f{d}{dt}\|\Dot{\Delta}_qu\|^2_{L^2}&+(\Dot{\Delta}_q(v\cdot\nabla u)|\Dot{\Delta}_qu)+\mu\|\nabla\Dot{\Delta}_qu\|^2_{L^2}+(\lambda+\mu)\|\dv\Dot{\Delta}_qu\|^2_{L^2}
+(\nabla\Dot{\Delta}_qc|\Dot{\Delta}_qu)\\&=(\Dot{\Delta}_qf|\Dot{\Delta}_qu),
\end{split}
\end{equation}
\begin{equation}
\label{14}\begin{split}
&\f{1}{2}\f{d}{dt}\|\Dot{\Delta}_qH\|^2_{L^2}+(\Dot{\Delta}_q(v\cdot\nabla H)|\Dot{\Delta}_qH)
+\|\Lambda\Dot{\Delta}_qH\|^2_{L^2}=(\Dot{\Delta}_qg|\Dot{\Delta}_qH),
\end{split}
\end{equation}
and
\begin{equation}
\label{15}\begin{split}
&\f{1}{2}\f{d}{dt}\|\nabla\Dot{\Delta}_qH\|^2_{L^2}+(\Dot{\Delta}_q(v\cdot\nabla H)|-\Delta\Dot{\Delta}_qH)
+\|\Delta\Dot{\Delta}_qH\|^2_{L^2}=(\Dot{\Delta}_qg|-\Delta\Dot{\Delta}_qH).
\end{split}
\end{equation}
Next, we derive an identity involving
$(\nabla\Dot{\Delta}_qc|\Dot{\Delta}_qu)$. For this purpose, we
apply the operator $\nabla$ to the first equation in
\eqref{equ:heatflow-new3} and take the $L^2$ scalar product with
$\Dot{\Delta}_qu$, then calculate the scalar product of the second
equation in \eqref{equ:heatflow-new3}  with $\nabla\Dot{\Delta}_qc$,
and then by summing up the results, we get
\begin{equation} \label{16}\begin{split}
&\f{d}{dt}(\nabla\Dot{\Delta}_qc|\Dot{\Delta}_qu)-\|\dv\Dot{\Delta}_qu\|^2_{L^2}+(\lambda+2\mu)(\dv\Dot{\Delta}_qu|\Delta\Dot{\Delta}_qc)
+\|\nabla\Dot{\Delta}_qc\|^2_{L^2}\\
&\qquad+(\nabla\Dot{\Delta}_q(v\cdot\nabla c)|\Dot{\Delta}_qu)+(\Dot{\Delta}_q(v\cdot\nabla u)|\nabla\Dot{\Delta}_qc)\\
&\quad=(\nabla\Dot{\Delta}_qf|\Dot{\Delta}_qu)+(\Dot{\Delta}_qg|\nabla\Dot{\Delta}_qc).
\end{split}
\end{equation}
We
now define
\begin{equation*}
\begin{split}
\alpha_q^2&=\|\Dot{\Delta}_qc\|^2_{L^2}+(\lambda+2\mu)A\|\nabla\Dot{\Delta}_qc\|^2_{L^2}+\|\Dot{\Delta}_qu\|^2_{L^2}+\|\Dot{\Delta}_qH\|^2_{L^2}+\|\nabla\Dot{\Delta}_qH\|^2_{L^2}
\\&\quad+2A(\nabla\Dot{\Delta}_qc|\Dot{\Delta}_qu),
\end{split}
\end{equation*}
where $A=\frac{\lambda+\mu}{2}>0$. Then, there exist two positive constants $c_1$ and $c_2$ such that
\begin{equation*}
\begin{split}
c_1\alpha_q^2\leq\|\Dot{\Delta}_qc\|^2_{L^2}+\|\nabla\Dot{\Delta}_qc\|^2_{L^2}+\|\Dot{\Delta}_qu\|^2_{L^2}+\|\Dot{\Delta}_qH\|^2_{L^2}+\|\nabla\Dot{\Delta}_qH\|^2_{L^2}
\leq c_2\alpha_q^2,
\end{split}
\end{equation*}
 for $M\in(1/(2\mu+\lambda),2/(\mu+\lambda))$,  we have
\begin{equation*}
\begin{split}
\big|2(\nabla\Dot{\Delta}_qc|\Dot{\Delta}_qu)\big|\leq M\|\Dot{\Delta}_qu\|^2_{L^2}+\frac{1}{M}\|\nabla\Dot{\Delta}_qc\|^2_{L^2}.
\end{split}
\end{equation*}
Thus, \begin{equation*}\label{equ:heatflow-new3}
\alpha_q\approx\left\{
\begin{aligned}{}
&\|(\Dot{\Delta}_qc,\Dot{\Delta}_qu,\Dot{\Delta}_qH)\|_{L^2}, \quad\quad \hbox{for} \quad q\leq q_0,\\
&\|(\nabla\Dot{\Delta}_qc,\Dot{\Delta}_qu,\nabla\Dot{\Delta}_qH)\|_{L^2}, \quad  \hbox{for} \quad q> q_0.
\end{aligned}
\right.
\end{equation*}
Combining with \eqref{11}-\eqref{16}, it yields, with the help of Proposition \ref{p29}, that
\begin{equation}\label{17}
\begin{split}
&\f{1}{2}\f{d}{dt}\alpha^2_q+(\mu+\lambda-A)\|\dv\Dot{\Delta}_qu\|^2_{L^2}+\mu\|\nabla\Dot{\Delta}_qu\|^2_{L^2}+A\|\nabla\Dot{\Delta}_q c\|^2_{L^2}
\\&\quad=-\big(\Dot{\Delta}_q(v\cdot\nabla c)|\Dot{\Delta}_qc\big)+(2\mu+\lambda)A\big(\Dot{\Delta}_q(v\cdot\nabla c)|\Delta\Dot{\Delta}_qc\big)-\big(\Dot{\Delta}_q(v\cdot\nabla u)|\Dot{\Delta}_qu\big)
\\&\qquad-\big(\Dot{\Delta}_q(v\cdot\nabla H)|\Dot{\Delta}_qH\big)
+\big(\Dot{\Delta}_q(v\cdot\nabla H)|\Delta\Dot{\Delta}_qH\big)-A\big(\nabla\Dot{\Delta}_q(v\cdot\nabla c)|\Dot{\Delta}_qu\big)
\\&\qquad-A\big(\Dot{\Delta}_q(v\cdot\nabla u)|\nabla\Dot{\Delta}_qc\big)
+\big(\Dot{\Delta}_qf|\Dot{\Delta}_qc\big)-(2\mu+\lambda)A\big(\Dot{\Delta}_qf|-\Delta\Dot{\Delta}_qc\big)+\big(\Dot{\Delta}_qf|\Dot{\Delta}_qu\big)
\\&\qquad+\big(\Dot{\Delta}_qg|\Dot{\Delta}_qH\big)
+\big(\Dot{\Delta}_qg|-\Delta\Dot{\Delta}_qH\big)+A\big(\nabla\Dot{\Delta}_qf|\Dot{\Delta}_qu\big)+A\big(\Dot{\Delta}_qg|\nabla\Dot{\Delta}_qc\big)
\\&\quad\lesssim \alpha_q\big(\|\Dot{\Delta}_qf\|_{L^2}+\|\nabla\Dot{\Delta}_qf\|_{L^2}+\|\Dot{\Delta}_qg\|_{L^2}+\|\Dot{\Delta}_qh\|_{L^2}+\|\nabla\Dot{\Delta}_qh\|_{L^2}
\\&\qquad+2^{-\frac{q}{2}}\gamma_qV'\|(c,u,H)\|_{\tilde{B}^{\f{1}{2},\f{3}{2}}_{2,1}\times
\dot{B}^{\f{1}{2}}_{2,1}\times \dot{B}^{\f{1}{2},\f{3}{2}}_{2,1}}\big).
\end{split}\end{equation}
Thus, it follows
\begin{align*}
  &\frac{1}{2}\frac{d}{dt}\alpha_q^2+c_0\min(2^{2q},1)\alpha_q^2\\
  &\quad\lesssim \gamma_q2^{-\frac{q}{2}}\Big[\|f\|_{\tilde{B}^{\f{1}{2},\f{3}{2}}_{2,1}} +\|g\|_{\dot{B}^{\f{1}{2}}_{2,1}}+\|h\|_{\dot{B}^{\f{1}{2},\f{3}{2}}_{2,1}}+V'\|(c,u,H)\|_{\tilde{B}^{\f{1}{2},\f{3}{2}}_{2,1}\times
\dot{B}^{\f{1}{2}}_{2,1}\times \dot{B}^{\f{1}{2},\f{3}{2}}_{2,1}}\Big)\Big]\alpha_q,
\end{align*}
which implies that
\begin{align*}
  &2^{\frac{q}{2}}\alpha_q+c_0\int_0^t\min(2^{2q},1)2^{\frac{q}{2}}\alpha_q(\tau)d\tau\\
  &\quad\lesssim2^{\frac{q}{2}}\alpha_q(0)+C\gamma_q\int_0^t\Big[\|f\|_{\tilde{B}^{\f{1}{2},\f{3}{2}}_{2,1}} +\|g\|_{\dot{B}^{\f{1}{2}}_{2,1}}+\|h\|_{\dot{B}^{\f{1}{2},\f{3}{2}}_{2,1}}+V'\sum_q2^{q(s-1)}\alpha_q\Big].
\end{align*}
Thus, by  Gronwall's inequality, we have
\begin{equation}\label{18}
\begin{split}
&\bigl\|c\bigr\|_{\widetilde{L}^{\infty}([0,t];\tilde{B}^{\f{1}{2},\f{3}{2}}_{2,1})}+\bigl\|u\bigr\|_{\widetilde{L}^{\infty}([0,t];
\dot{B}^{\f{1}{2}}_{2,1})}+\bigl\|H\bigr\|_{\widetilde{L}^{\infty}([0,t];\tilde{B}^{\f{1}{2},\f{3}{2}}_{2,1})}
\\&\qquad+\int_0^t\big\|c(\tau)\big\|_{\tilde{B}^{\f{5}{2},\f{3}{2}}_{2,1}}d\tau+\int_0^t\big\|u(\tau)\big\|_{
\tilde{B}^{\f{5}{2},\f{1}{2}}_{2,1}}d\tau+\int_0^t\big\|H(\tau)\big\|_{\tilde{B}^{\f{5}{2},\f{3}{2}}_{2,1}}d\tau
\\&\quad \lesssim e^{CV(t)}\Big( \|c_{0}\|_{\tilde{B}^{\f{1}{2},\f{3}{2}}_{2,1}}+\|u_{0}\|_{
\dot{B}^{\f{1}{2}}_{2,1}}+\|H_{0}\|_{\tilde{B}^{\f{1}{2},\f{3}{2}}_{2,1}}+\int_0^t\big\|f(\tau)\big\|_{\tilde{B}^{\f{1}{2},\f{3}{2}}_{2,1}}d\tau
\\&\qquad+\int_0^t\big\|g(\tau)\big\|_{\dot{B}^{\f{1}{2}}_{2,1}}d\tau+\int_0^t\big\|h(\tau)\big\|_{\tilde{B}^{\f{1}{2},\f{3}{2}}_{2,1}}d\tau\Big).
\end{split}
\end{equation}
Based on the damping effect for $c$, we   get the smoothing
effect of $u$  by considering
\eqref{equ:heatflow-new3} with $\nabla c$ being seen as a source
term. Furthermore, we also exploit the smoothing
effect of $H$ by heat equation of $\eqref{equ:heatflow-new3}_{3}$. Thanks to \eqref{18}, it suffices to state the proof
for the high frequencies only for  $u$ and $H$. From \eqref{13} and Proposition \ref{p29}, we have
\begin{equation*}
\label{19}\begin{split}
&\f{1}{2}\f{d}{dt}\|\Dot{\Delta}_qu\|^2_{L^2}+C2^{2q}\|\Dot{\Delta}_qu\|^2_{L^2}
\\&\quad=-(\nabla\Dot{\Delta}_qc|\Dot{\Delta}_qu)-(\Dot{\Delta}_q(v\cdot\nabla u)|\Dot{\Delta}_qu)+(\Dot{\Delta}_qf|\Dot{\Delta}_qu)
\\&\quad\lesssim \|\Dot{\Delta}_qu\|_{L^2}\Big(2^{q}\|\Dot{\Delta}_qc\|_{L^{2}}+\|\Dot{\Delta}_qg\|_{L^{2}}+ V'(t)\gamma_{q}2^{-\frac{q}{2}}\|u\|_{
\dot{B}^{\frac{1}{2}}_{2,1}}\Big).
\end{split}
\end{equation*}
It follows that
\begin{equation*}
\label{19}\begin{split}
&\f{d}{dt}\sum_{q\geq q_0}2^{\frac{q}{2}}\|\Dot{\Delta}_qu\|_{L^2}+C\sum_{q\geq q_0}2^{\frac{q}{2}}2^{2q}\|\Dot{\Delta}_qu\|_{L^2}
\\&\quad\lesssim \sum_{q\geq q_0}2^{\frac{q}{2}}\Big(2^{q}\|\Dot{\Delta}_qc\|_{L^{2}}+\|\Dot{\Delta}_qg\|_{L^{2}}+ V'(t)\gamma_{q}2^{-\frac{q}{2}}\|u\|_{
\dot{B}^{\frac{1}{2}}_{2,1}}\Big)
\\&\quad\lesssim \sum_{q\geq q_0}2^{\frac{q}{2}}2^{q}\|\Dot{\Delta}_qc\|_{L^{2}}+\|\Dot{\Delta}_qg\|_{\dot{B}^{\frac{1}{2}}_{2,1}}+ V'(t)\|u\|_{
\dot{B}^{\frac{1}{2}}_{2,1}},
\end{split}
\end{equation*}
which implies, with the help of \eqref{18}, that
\begin{equation}
\label{20}\begin{split}
&\int_0^t\sum_{q\geq q_{0}}2^{\frac{q}{2}}2^{2q}\|\Dot{\Delta}_qu(\tau)\|_{L^2}d\tau
\\&\quad\lesssim \|u_0\|_{
\dot{B}^{\frac{q}{2}}_{2,1}}+\int_0^t\sum_{q\geq q_0}2^{\frac{q}{2}}2^{q}\|\Dot{\Delta}_qc(\tau)\|_{L^{2}}d\tau+\int_0^t\|g(\tau)\|_{\dot{B}^{\frac{1}{2}}_{2,1}}d\tau+ V(t)\sup_{\tau\in[0,t]}\|u\|_{
\dot{B}^{\frac{1}{2}}_{2,1}}
\\&\quad\lesssim e^{CV(t)}\Big( \|c_{0}\|_{\tilde{B}^{\f{1}{2},\f{3}{2}}_{2,1}}+\|u_{0}\|_{
\dot{B}^{\f{1}{2}}_{2,1}}+\|H_{0}\|_{\tilde{B}^{\f{1}{2},\f{3}{2}}_{2,1}}+\int_0^t\big\|f(\tau)\big\|_{\tilde{B}^{\f{1}{2},\f{3}{2}}_{2,1}}d\tau
\\&\qquad+\int_0^t\big\|g(\tau)\big\|_{\dot{B}^{\f{1}{2}}_{2,1}}d\tau+\int_0^t\big\|h(\tau)\big\|_{\tilde{B}^{\f{1}{2},\f{3}{2}}_{2,1}}d\tau\Big).
\end{split}
\end{equation}
Similarly, from \eqref{14},  Proposition \ref{p29} and Remark  \ref{2.5}, we have
\begin{equation}
\label{21}\begin{split}
&\int_0^t\sum_{q\geq q_{0}}2^{\frac{3q}{2}}2^{2q}\|\Dot{\Delta}_qH(\tau)\|_{L^2}d\tau
\\&\quad\lesssim \|H_0\|_{
\dot{B}^{\frac{3}{2}}_{2,1}}+\int_0^t\|g(\tau)\|_{\dot{B}^{\frac{3}{2}}_{2,1}}d\tau+ V(t)\sup_{\tau\in[0,t]}\|H\|_{
\dot{B}^{\frac{1}{2},\frac{3}{2}}_{2,1}}
\\&\quad\lesssim \|H_0\|_{
\dot{B}^{\frac{1}{2},\frac{3}{2}}_{2,1}}+\int_0^t\|h(\tau)\|_{\dot{B}^{\frac{1}{2},\frac{3}{2}}_{2,1}}d\tau+ V(t)\sup_{\tau\in[0,t]}\|H\|_{
\dot{B}^{\frac{1}{2},\frac{3}{2}}_{2,1}}
\\&\quad\lesssim e^{CV(t)}\Big( \|c_{0}\|_{\tilde{B}^{\f{1}{2},\f{3}{2}}_{2,1}}+\|u_{0}\|_{
\dot{B}^{\f{1}{2}}_{2,1}}+\|H_{0}\|_{\tilde{B}^{\f{1}{2},\f{3}{2}}_{2,1}}+\int_0^t\big\|f(\tau)\big\|_{\tilde{B}^{\f{1}{2},\f{3}{2}}_{2,1}}d\tau
\\&\qquad+\int_0^t\big\|g(\tau)\big\|_{\dot{B}^{\f{1}{2}}_{2,1}}d\tau+\int_0^t\big\|h(\tau)\big\|_{\tilde{B}^{\f{1}{2},\f{3}{2}}_{2,1}}d\tau\Big).
\end{split}
\end{equation}
Combining with \eqref{20}-\eqref{21} and \eqref{18}, we finally conclude that \eqref{4.21}. Thus,  we complete the
proof of Proposition \ref{prop3.1}.

\section{Global existence for initial data near equilibrium}
 \ \ \  In this section, we are going to show that if the initial
 data satisfy
$$\|c_0\|_{\tilde{B}^{\f{1}{2},\f{3}{2}}_{2,1}}+\|u_0\|_{\dot{B}^{\f{1}{2}}_{2,1}}+\|H_0\|_{\dot{B}^{\f{1}{2},\frac{3}{2}}_{2,1}}\le
\eta$$ for some sufficiently small $\eta$,  then there exists a
positive constant $M$ such that
$$\|(c, u, H)\|_{\mathcal{E}^{\f{3}{2}}}\le M\eta,$$
where \begin{equation*}\label{eq:conver2}
\begin{split} \mathcal{E}^{\f{3}{2}}\eqdefa
&\big(L^1(\mathbb{R}^{+}; \tilde{B}^{\frac{5}{2},\frac{3}{2}}_{2,1})\cap
C(\mathbb{R}^{+};\tilde{B}^{\frac{1}{2},\frac{3}{2}}_{2,1})\big)\times\big(
L^1(\mathbb{R}^{+}; \dot{B}^{\frac{5}{2}}_{2,1})\cap
C(\mathbb{R}^{+};\dot{B}^{\frac{1}{2}}_{2,1})\big)\\
&\quad\times\big(L^1(\mathbb{R}^{+}; \tilde{B}^{\frac{5}{2},\frac{7}{2}}_{2,1})\cap
C(\mathbb{R}^{+};\tilde{B}^{\frac{1}{2},\frac{3}{2}}_{2,1})\big).
\end{split}
\end{equation*}
This uniform estimate will enable us to
 extend the local solution $(c, u, H)$ obtained within an iterative scheme as in \cite{Dan1} to  a global one. To this
end, we use a contradiction argument. Define
$$T_0=\sup\big\{T\in[0,\infty): \|(c,u,H,\theta)\|_{\mathcal{E}_T^{\f{3}{2}}}\le
M\eta\big\},$$ with $M$ to be determined later. Suppose
that $T_0<\infty$. We apply  linear estimates in Proposition \ref{prop3.1}
to the solutions of the reformulated system
\eqref{equ:heatflow-new} such that for all $t\in [0,T_0]$, the
following estimates hold
\begin{equation}\label{41}
\begin{split}
&\bigl\|c\bigr\|_{\widetilde{L}^{\infty}([0,t];\tilde{B}^{\f{1}{2},\f{3}{2}}_{2,1})}+\bigl\|u\bigr\|_{\widetilde{L}^{\infty}([0,t];
\dot{B}^{\f{1}{2}}_{2,1})}+\bigl\|H\bigr\|_{\widetilde{L}^{\infty}([0,t];\tilde{B}^{\f{1}{2},\f{3}{2}}_{2,1})}
\\&\qquad+\int_0^t\big\|c(\tau)\big\|_{\tilde{B}^{\f{5}{2},\f{3}{2}}_{2,1}}d\tau+\int_0^t\big\|u(\tau)\big\|_{
\dot{B}^{\f{5}{2}}_{2,1}}d\tau+\int_0^t\big\|H(\tau)\big\|_{\tilde{B}^{\f{5}{2},\f{7}{2}}_{2,1}}d\tau
\\&\quad \lesssim e^{CV(t)}\Big( \|c_{0}\|_{\tilde{B}^{\f{1}{2},\f{3}{2}}_{2,1}}+\|u_{0}\|_{
\dot{B}^{\f{1}{2}}_{2,1}}+\|H_{0}\|_{\tilde{B}^{\f{1}{2},\f{3}{2}}_{2,1}}+\int_0^t\big\|f(\tau)\big\|_{\tilde{B}^{\f{1}{2},\f{3}{2}}_{2,1}}d\tau
\\&\qquad+\int_0^t\big\|g(\tau)\big\|_{\dot{B}^{\f{1}{2}}_{2,1}}d\tau+\int_0^t\big\|h(\tau)\big\|_{\tilde{B}^{\f{1}{2},\f{3}{2}}_{2,1}}d\tau\Big).
\end{split}
\end{equation} where \begin{equation*}
V(T_0)=\int_0^{T_0}\|u(\tau)\|_{\dot{B}^{\f{5}{2}}_{2,1}}d\tau.
\label{VV}
\end{equation*}
In what follows, we derive some estimates for the nonlinear terms
$f,g$ and $h$. First, by Proposition \ref{p25},  we have
\begin{equation}\label{42}
\begin{split}\|f\|_{L^1_{T_0}(\tilde{B}^{\f{1}{2},\f{3}{2}}_{2,1})}&\lesssim\|c\|_{L^\infty_{T_0}(\tilde{B}^{\frac{1}{2},\f 32}_{2,1})}\|\Dv u\|_{L^1_{T_0}(\dot{B}^{\f 32}_{2,1})}\\
&\lesssim M^2\eta^2.
\end{split}
\end{equation}
Next, we bound the term $g$. By the embedding $\tilde{B}^{\frac{1}{2},\f
32}_{2,1}\hookrightarrow \dot{B}^\frac 32_{2,1}\hookrightarrow
L^\infty$ and Proposition \ref{p27}, we get
$$\|L_1(c)\|_{\dot{B}^\frac 32_{2,1}}\leq C_0(\|c\|_{L^\infty})\|c\|_{\dot{B}^\frac 3 2_{2,1}}\lesssim \|c\|_{\tilde{B}^{\frac{1}{2},\f 32}_{2,1}},$$
thus
\begin{equation*}\label{318}
\begin{split}
\|L_1(c)\mathcal{A}u\|_{L^1_{T_0}(\dot{B}^{\f{1}{2}}_{2,1})}
&\lesssim\|L_1(c)\|_{L^\infty_{T_0}(\dot{B}^\frac 32_{2,1})}
\|\mathcal{A}u\|_{L^1_{T_0}(\dot{B}^{\f{1}{2}}_{2,1})}\\
&\lesssim\|c\|_{L^\infty_{T_0}(\dot{B}^\frac 32_{2,1})}
\|\mathcal{A}u\|_{L^1_{T_0}(\dot{B}^{\f{1}{2}}_{2,1})}\\
&\lesssim\|c\|_{L^\infty_{T_0}(\tilde{B}^{\frac{1}{2},\f 32}_{2,1})}\|u\|_{L^1_{T_0}(\dot{B}^{\f{5}{2}}_{2,1})}\\
&\lesssim M^2\eta^2.
\end{split}
\end{equation*}
Similarly,
\begin{equation*}\label{318}
\begin{split}&\Big\|L_3(c)\big(\frac{1}{2}\nabla |H|^{2}-H\cdot\nabla H\big)\Big\|_{L^1_{T_0}(\dot{B}^{\f{1}{2}}_{2,1})}\\
&\quad\lesssim\big(1+\|L_1(c)\|_{L^\infty_{T_0}(\dot{B}^\frac 32_{2,1})}\big)
\|H\cdot\nabla H\|_{L^1_{T_0}(\dot{B}^{\f{1}{2}}_{2,1})}\\
&\quad\lesssim\big(1+\|L_1(c)\|_{L^\infty_{T_0}(\dot{B}^\frac
32_{2,1})}\big)\|H\|_{L^\infty_{T_0}(\dot{B}^{\frac
{1}{2}}_{2,1})}
\|\nabla H\|_{L^1_{T_0}(\dot{B}^{\f{3}{2}}_{2,1})}\\
&\quad\lesssim\big(1+\|L_1(c)\|_{L^\infty_{T_0}(\dot{B}^\frac
32_{2,1})}\big)\|H\|_{L^\infty_{T_0}(\dot{B}^{\frac
{1}{2}}_{2,1})} \| H\|_{L^1_{T_0}(\dot{B}^{\f{5}{2}}_{2,1})}
\\&\quad\le C\big(1+\|L_1(c)\|_{L^\infty_{T_0}(\dot{B}^{\frac
32}_{2,1})}\big)\|H\|_{L^\infty_{T_0}(\dot{B}^{\frac
{1}{2},\frac
{3}{2}}_{2,1})} \| H\|_{L^1_{T_0}(\dot{B}^{\f{5}{2},\f{7}{2}}_{2,1})}
\\&\quad\lesssim(1+M\eta)M^2\eta^2.
\end{split}
\end{equation*}
From  Propositions
\ref{p26}-\ref{p27} and Remark \ref{2.5}, we conclude that
\begin{equation*}\label{43}
\begin{split}
\|L_2(c)\na c\|_{L^1_{T_0}(\dot{B}^{\f{1}{2}}_{2,1})}&\lesssim
\|L_2(c)\|_{L^2_{T_0}(\tilde{B}^{\frac{3}{2}}_{2,1})} \|
c\|_{L^2_{T_0}(\dot{B}^{\f{3}{2}}_{2,1})}
\\&\lesssim\|c\|^{2}_{L^2_{T_0}(\dot{B}^{\f{3}{2}}_{2,1})}
\\&\lesssim\|c\|_{L^\infty_{T_0}(\dot{B}^{\f{1}{2},\frac{3}{2}}_{2,1})}\|c\|_{L^1_{T_0}(\dot{B}^{\f{5}{2},\frac{3}{2}}_{2,1})}
\\&\lesssim M^2\eta^2.
\end{split}
\end{equation*}
Therefore,
\begin{equation}\label{43}
\|g\|_{L^1_{T_0}(\dot{B}^{\f{1}{2}}_{2,1})}\le
C(1+M\eta)M^2\eta^2.
\end{equation}
Finally,  we  bound the term $h$ as follows. Employing  Proposition
\ref{p26}, we infer that
\begin{equation*}\label{44443}
\begin{split}
\big\|H\cdot\nabla u-\dv u
H\big\|_{L^1_{T_0}(\dot{B}^{\f{1}{2},\frac{3}{2}}_{2,1})}&\lesssim\int_0^{T_0}\|H\|_{\dot{B}^{\f{1}{2},\frac{3}{2}}_{2,1}}\|\nabla u\|_{\dot{B}^{\f{3}{2}}_{2,1}}d\tau\\
&\lesssim\|H\|_{L^\infty_{T_0}(\dot{B}^{\f{1}{2},\frac{3}{2}}_{2,1})}\|u\|_{L^1_{T_0}(\dot{B}^{\f{5}{2}}_{2,1})}\\
&\lesssim M^2\eta^2,
\end{split}
\end{equation*}
\begin{equation*}\label{44443}
\begin{split}
&\Big\|\nabla\times \big(L_3(c)(\nabla\times H)\times H\big)\Big\|_{L^1_{T_0}(\dot{B}^{\f{1}{2},\frac{3}{2}}_{2,1})}\\
&\quad\lesssim
\Big\|\big(L_3(c)(\nabla\times H)\times H\big)\Big\|_{L^1_{T_0}(\dot{B}^{\f{3}{2},\frac{5}{2}}_{2,1})}\\
&\quad\lesssim\big(1+\|L_1(c)\|_{L^\infty_{T_0}(\dot{B}^{\frac
32}_{2,1})}\big)\big\|(\nabla\times H)\times H\big\|_{L^1_{T_0}(\dot{B}^{\f{3}{2},\frac{5}{2}}_{2,1})}
\\
&\quad\lesssim\big(1+\|L_1(c)\|_{L^\infty_{T_0}(\dot{B}^{\frac
32}_{2,1})}\big)\| H\|_{L^\infty_{T_0}(\dot{B}^{\f{3}{2},\frac{3}{2}}_{2,1})}\|\nabla\times H\|_{L^1_{T_0}(\dot{B}^{\f{3}{2},\frac{5}{2}}_{2,1})}
\\
&\quad\lesssim\big(1+\|L_1(c)\|_{L^\infty_{T_0}(\dot{B}^{\frac
32}_{2,1})}\big)\| H\|_{L^\infty_{T_0}(\dot{B}^{\f{1}{2},\frac{3}{2}}_{2,1})}\| H\|_{L^1_{T_0}(\dot{B}^{\f{5}{2},\frac{7}{2}}_{2,1})}\\
&\quad\lesssim(1+M\eta)M^2\eta^2.
\end{split}
\end{equation*}
Hence, we gather that
\begin{equation}\label{44}
\|h\|_{L^1_{T_0}(\dot{B}^{\f{1}{2}}_{2,1})}\lesssim(1+M\eta)M^2\eta^2.
\end{equation}
 Substituting
\eqref{42}-\eqref{44} into  \eqref{41} , we obtain that
\begin{equation}\label{322}
\begin{split}
\|(c,u,H)\|_{\mathcal{E}_{T_0}^{\f{3}{2}}}\leq
C_1e^{C_1M\eta}\left(1+M\eta\right)M^2\eta^2.
\end{split}
\end{equation}
Choose $M=8C_1$, for $\eta$ small enough such that
\begin{equation}
e^{C_1M \eta}\leq 2, \quad
\left(1+M\eta\right)M^2\eta\leq2,
\end{equation}
which implies that
$$\|(c,u,H)\|_{\mathcal{E}_{T_{0}}^{\f{3}{2}}}\le \frac{1}{2}M\eta. $$ This is a
contradiction with the definition of $T_0$. As a consequence, we
conclude that $T_0=\infty$. Based on the above uniform estimates,   employing a classical Friedrich's approximation and compactness method
(cf. \cite{Dan2,Dan1,Dan3}), we can establish the  global existence
of  strong solutions of the system \eqref{equ:heatflow-new}. Here,
we omit it. This completes the proof of the existence in Theorem \ref{th:main1}.
\section{Uniqueness}
 \ \  In this section, we will address  uniqueness of  strong solutions to the system \eqref{equ:heatflow-new}. For this purpose, suppose that $(c_i,u_i,
H_i)_{i=1,2}$ in  $\mathcal{E}_T^{\f{3}{2}}$  solve \eqref{equ:heatflow-new} with
the same initial data.  Define $$\delta c=c_2-c_1, \ \ \delta u=u_2-u_1,\ \  \delta
H=H_2-H_1.$$ Then, $(\delta
c,\delta u, \delta H)$ satisfies the following system
\begin{equation}\label{51}
\left\{
\begin{aligned}{}
&\p_t \delta c+u_2\cdot\na \delta c+\dv \delta u=\delta \mathrm{f},\\
&\p_t \delta u+u_2\cdot\na  \delta u-\mathcal{A}\delta u+\nabla  \delta c=\delta \mathrm{g},\\
&\p_t \delta H+u_2\cdot\nabla \delta H- \Delta \delta H=\delta \mathrm{h},
\end{aligned}
\right.
\end{equation}
where
\begin{equation*}\label{eq:conver2}
\delta \mathrm{f}=-\delta u\cdot\nabla c_1-\delta c\Dv u_2-c_1\Dv\delta u,\quad
\delta \mathrm{g}=\sum_{i=1}^{8}\delta \mathrm{g}_{i},\quad
\delta \mathrm{h}=\sum_{i=1}^{8}\delta \mathrm{h}_{i},
\end{equation*}
with
\begin{align*}
\delta \mathrm{g}_{1}&= - u_2\cdot\nabla \delta u, & \qquad  \delta \mathrm{h}_{1}&=\delta u\cdot\nabla H_{1}, \\
\delta \mathrm{g}_{2}&=-\delta u\cdot\nabla u_1, & \qquad \delta \mathrm{h}_{2}&=H_2\nabla\delta u,\\
\delta \mathrm{g}_{3}&=-\big(L_{1}(c_2)-L_{1}(c_1)\big)\mathcal{A} u_2, & \qquad \delta \mathrm{h}_{3}&=\delta H\cdot\nabla u_1,\\
\delta \mathrm{g}_{4}&=-L_{1}(c_1)\mathcal{A} \delta u, & \qquad  \delta \mathrm{h}_{4}&=-H_2\Dv\delta u,\\
\delta \mathrm{g}_{5}&=-\nabla \big(K_0(c_2)-K_0(c_1)\big),& \qquad  \delta \mathrm{h}_{5}&=-\delta H\Dv u_1, \\
\delta \mathrm{g}_{6}&=-\big(L_{3}(c_2)-L_{3}(c_1)\big)\nabla H_2\cdot H_2,& \qquad \delta \mathrm{h}_{6}&=-\nabla\times\Big[\big(L_3(c_2)-L_3(c_1)\big)(\nabla\times H_2)\times H_2\Big],\\
\delta \mathrm{g}_{7}&=-L_{3}(c_1)\nabla\delta H\cdot H_2, & \qquad \delta \mathrm{h}_{7}&=\nabla\times\big[L_3(c_1)(\nabla\times \delta H)\times H_2\big],\\
\delta \mathrm{g}_{8}&=-L_{3}(c_1)\nabla H_1\cdot \delta H, & \qquad \delta \mathrm{h}_{8}&=\nabla\times\big[L_3(c_1)(\nabla\times H_1)\times \delta H\big],
\end{align*}
 and
$$ K_0(z)=\int_{0}^{z}L_{2}(y)dy.$$
Applying Proposition \ref{prop3.1} to the system \eqref{51}, we get
\begin{equation}\label{52}
\begin{split}
&\|(\delta c, \delta u, \delta H)\|_{\mathcal{E}_{T}^{\f{1}{2}}}
\\&\quad \lesssim e^{\int_0^t\|u_2(\tau)\|_{\dot{B}^{\f{5}{2}}_{2,1}}d\tau}\Big( \big\|\delta \mathrm{f}\big\|_{L_{T}^{1}(\tilde{B}^{-\f{1}{2},\f{1}{2}}_{2,1})}
+\big\|\delta \mathrm{g}\big\|_{L_{T}^{1}(\dot{B}^{-\f{1}{2}}_{2,1})}+\big\|\delta \mathrm{h}\big\|_{L_{T}^{1}(\tilde{B}^{-\f{1}{2},\f{1}{2}}_{2,1})}\Big).
\end{split}
\end{equation}
Let us observe that $\partial_tc_i\in L^1_{loc}(\dot{B}_{2,1}^{\frac{1}{2}})$, and hence $c_i\in C(\dot{B}_{2,1}^{\frac{1}{2}})\cap L^\infty(\dot{B}_{2,1}^{\frac{3}{2}})(i=1,2)$. This entails $c_i\in C([0,\infty)\times \mathbb{R}^3)$. On the other
hand, if $\eta$ is sufficiently small, we have
$$|c_1(t,x)|\le \frac{1}{4}
\quad\textrm{for all}\quad t\ge 0\textrm{ and }x\in\mathbb{R}^3.$$
Continuity in time for $c_2$ thus yields the existence of a time $T>0$ such that
$$\|c_i(t)\|_{L^\infty}\le \frac{1}{2} \quad\textrm{for}\quad i=1,2\textrm{  and  }t\in[0,T].$$
Employing Propositions
\ref{p26}-\ref{p27} and Remarks \ref{2.3}-\ref{2.5}, we easily infer that
\begin{equation*}
 \begin{split}
&\|\delta\mathrm{f}\|_{L_T^1(\tilde{B}_{2,1}^{-\f{1}{2},\f{1}{2}})}\\&\quad\lesssim \|c_1\|_{L^\infty_T(\tilde{B}_{2,1}^{\f12,\f32})}\|\delta u\|_{L^1_T(\dot{B}_{2,1}^{\f32})}+\|\Dv u_2\|_{L^1_T(\dot{B}_{2,1}^{\f32})}\|\delta c\|_{L^\infty_T(\tilde{B}_{2,1}^{-\f{1}{2},\f{1}{2}})},
 \end{split}
\end{equation*}
\begin{equation*}
 \begin{split}
&\|\delta\mathrm{g}_1\|_{L_T^1(\dot{B}_{2,1}^{-\f{1}{2}})}+\|\delta\mathrm{g}_2\|_{L_T^1(\dot{B}_{2,1}^{-\f{1}{2}})}
+\|\delta\mathrm{g}_3\|_{L_T^1(\dot{B}_{2,1}^{-\f{1}{2}})}+
\|\delta\mathrm{g}_4\|_{L_T^1(\dot{B}_{2,1}^{-\f{1}{2}})}\\&\quad\lesssim \|u_2\|_{L^2_T(\dot{B}_{2,1}^{\f32})}\|\nabla\delta u\|_{L^2_T(\dot{B}_{2,1}^{-\f12})}+\|\delta u\|_{L^2_T(\dot{B}_{2,1}^{\f12})}\|\nabla u_1\|_{L^2_T(\dot{B}_{2,1}^{\f12})}
\\&\qquad+\Big(1+\|c_1\|_{L^\infty_T(\dot{B}_{2,1}^{\f32})}+\|c_2\|_{L^\infty_T(\dot{B}_{2,1}^{\f32})}\Big)\|\delta c\|_{L^\infty_T(\dot{B}_{2,1}^{\f12})}\|\nabla^2u_2\|_{L^1_T(\dot{B}_{2,1}^{\f12})}
\\&\qquad+T\Big(\|c_1\|_{L^\infty_T(\dot{B}_{2,1}^{\f32})}+\|c_2\|_{L^\infty_T(\dot{B}_{2,1}^{\f32})}\Big)\|\delta c\|_{L^\infty_T(\dot{B}_{2,1}^{\f12})}
\\&\qquad+\| c_1\|_{L^\infty_T(\dot{B}_{2,1}^{\f32})}\|\nabla^2\delta u\|_{L^1_T(\dot{B}_{2,1}^{-\f12})}
\\&\quad\lesssim \|u_2\|_{L^2_T(\dot{B}_{2,1}^{\f32})}\|\delta u\|_{L^2_T(\dot{B}_{2,1}^{\f12})}+\|\delta u\|_{L^2_T(\dot{B}_{2,1}^{\f12})}\| u_1\|_{L^2_T(\dot{B}_{2,1}^{\f32})}
\\&\qquad+\Big(1+\|c_1\|_{L^\infty_T(\dot{B}_{2,1}^{\f32})}+\|c_2\|_{L^\infty_T(\dot{B}_{2,1}^{\f32})}\Big)\|\delta c\|_{L^\infty_T(\tilde{B}_{2,1}^{-\f12,\f12})}\|u_2\|_{L^1_T(\dot{B}_{2,1}^{\f52})}
\\&\qquad+T\Big(\|c_1\|_{L^\infty_T(\dot{B}_{2,1}^{\f32})}+\|c_2\|_{L^\infty_T(\dot{B}_{2,1}^{\f32})}\Big)\|\delta c\|_{L^\infty_T(\tilde{B}_{2,1}^{-\f12,\f12})}
\\&\qquad+\| c_1\|_{L^\infty_T(\tilde{B}_{2,1}^{\f12,\f32})}\|\delta u\|_{L^1_T(\dot{B}_{2,1}^{\f32})},
 \end{split}
\end{equation*}
\begin{equation*}
 \begin{split}
&\|\delta\mathrm{g}_5\|_{L_T^1(\dot{B}_{2,1}^{-\f{1}{2}})}+\|\delta\mathrm{g}_6\|_{L_T^1(\dot{B}_{2,1}^{-\f{1}{2}})}
+\|\delta\mathrm{g}_7\|_{L_T^1(\dot{B}_{2,1}^{-\f{1}{2}})}+
\|\delta\mathrm{g}_8\|_{L_T^1(\dot{B}_{2,1}^{-\f{1}{2}})}\\&\quad\lesssim \Big(1+\|c_1\|_{L^\infty_T(\dot{B}_{2,1}^{\f32})}+\|c_2\|_{L^\infty_T(\dot{B}_{2,1}^{\f32})}\Big)\|\delta c\|_{L^\infty_T(\dot{B}_{2,1}^{\f12})}\|\nabla H_{2}\cdot H_2\|_{L^1_T(\dot{B}_{2,1}^{\f12})}\\
&\qquad+\Big(1+\|c_1\|_{L^\infty_T(\dot{B}_{2,1}^{\f32})}\Big)\|\nabla H_{2}\cdot \delta H\|_{L^1_T(\dot{B}_{2,1}^{-\f12})}\\
&\qquad+\Big(1+\|c_1\|_{L^\infty_T(\dot{B}_{2,1}^{\f32})}\Big)\|H_2\cdot\nabla\delta H\|_{L^1_T(\dot{B}_{2,1}^{-\f12})}
\\&\quad\lesssim\Big(1+\|c_1\|_{L^\infty_T(\dot{B}_{2,1}^{\f32})}+\|c_2\|_{L^\infty_T(\dot{B}_{2,1}^{\f32})}\Big)\|\delta c\|_{L^\infty_T(\tilde{B}_{2,1}^{-\f12,\f12})}\| H_2\|_{L^\infty_T(\dot{B}_{2,1}^{\f12,\f32})}\| H_{2}\|_{L^1_T(\tilde{B}_{2,1}^{\f52,\f72})}\\
&\qquad+\Big(1+\|c_1\|_{L^\infty_T(\dot{B}_{2,1}^{\f32})}\Big)\| H_{2}\|_{L^1_T(\tilde{B}_{2,1}^{\f52,\f72})}\| \delta H\|_{L^\infty_T(\tilde{B}_{2,1}^{-\f12,\f12})}\\
&\qquad+\Big(1+\|c_1\|_{L^\infty_T(\dot{B}_{2,1}^{\f32})}\Big)\|H_2\|_{L^2_T(\tilde{B}_{2,1}^{\f32,\f52})}\|\delta H\|_{L^2_T(\tilde{B}_{2,1}^{\f12,\f32})},
 \end{split}
\end{equation*}
\begin{equation*}
 \begin{split}
&\|\delta\mathrm{h}_1\|_{L_T^1(\dot{B}_{2,1}^{-\f{1}{2}})}+\|\delta\mathrm{h}_2\|_{L_T^1(\dot{B}_{2,1}^{-\f{1}{2}})}
+\|\delta\mathrm{h}_3\|_{L_T^1(\dot{B}_{2,1}^{-\f{1}{2}})}+
\|\delta\mathrm{h}_4\|_{L_T^1(\dot{B}_{2,1}^{-\f{1}{2}})}+
\|\delta\mathrm{h}_5\|_{L_T^1(\dot{B}_{2,1}^{-\f{1}{2}})}\\
&\quad\lesssim \|\delta u\|_{L^\infty_T(\tilde{B}_{2,1}^{-\f12,-\f12})}\|\nabla H_1\|_{L^1_T(\tilde{B}_{2,1}^{\f32,\f52})}+\|\nabla\delta u\|_{L^2_T(\tilde{B}_{2,1}^{-\f12,-\f12})}\| H_2\|_{L^2_T(\tilde{B}_{2,1}^{\f32,\f52})}\\
&\qquad+\|\delta H\|_{L^\infty_T(\tilde{B}_{2,1}^{-\f12,\f12})}\|\nabla u_1\|_{L^1_T(\tilde{B}_{2,1}^{\f32,\f32})}\\
&\quad\lesssim \|\delta u\|_{L^\infty_T(\dot{B}_{2,1}^{-\f12})}\| H_1\|_{L^1_T(\tilde{B}_{2,1}^{\f52,\f72})}+\|\delta u\|_{L^2_T(\dot{B}_{2,1}^{\f12})}\| H_2\|_{L^2_T(\tilde{B}_{2,1}^{\f32,\f52})}\\
&\qquad+\|\delta H\|_{L^\infty_T(\tilde{B}_{2,1}^{-\f12,\f12})}\| u_1\|_{L^1_T(\dot{B}_{2,1}^{\f52})}\\
&\quad\lesssim \|\delta u\|_{L^\infty_T(\dot{B}_{2,1}^{-\f12})}\| H_1\|_{L^1_T(\tilde{B}_{2,1}^{\f52,\f72})}+\|\delta u\|_{L^2_T(\dot{B}_{2,1}^{\f12})}\| H_2\|_{L^2_T(\tilde{B}_{2,1}^{\f32,\f52})}\\
&\qquad+\|\delta H\|_{L^\infty_T(\tilde{B}_{2,1}^{-\f12,\f12})}\| u_1\|_{L^1_T(\dot{B}_{2,1}^{\f52})},
 \end{split}
\end{equation*}
and
\begin{equation*}
 \begin{split}
&\|\delta\mathrm{h}_6\|_{L_T^1(\dot{B}_{2,1}^{-\f{1}{2}})}+
\|\delta\mathrm{h}_7\|_{L_T^1(\dot{B}_{2,1}^{-\f{1}{2}})}+
\|\delta\mathrm{h}_8\|_{L_T^1(\dot{B}_{2,1}^{-\f{1}{2}})}\\&\quad\lesssim \Big\|\big(L_3(c_2)-L_3(c_1)\big)(\nabla\times H_2)\times H_2\Big\|_{L^1_T(\tilde{B}_{2,1}^{\f12,\f32})}\\
&\qquad+\Big\|L_3(c_1)(\nabla\times \delta H)\times H_2\Big\|_{L^1_T(\tilde{B}_{2,1}^{\f12,\f32})}\\
&\qquad+\Big\|L_3(c_1)(\nabla\times H_1)\times \delta H\Big\|_{L^1_T(\tilde{B}_{2,1}^{\f12,\f32})}\\
&\quad\lesssim \Big(1+\|c_1\|_{L^\infty_T(\dot{B}_{2,1}^{\f32}))}+\|c_2\|_{L^\infty_T(\dot{B}_{2,1}^{\f32}))}\Big)\|\delta c\|_{L^\infty_T(\tilde{B}_{2,1}^{\f12,\f12})}\big\|(\nabla\times H_2)\times H_2\big\|_{L^1_T(\tilde{B}_{2,1}^{\f32,\f52})}\\
&\qquad+\Big(1+\|c_1\|_{L^\infty_T(\dot{B}_{2,1}^{\f32})}\Big)\big\|(\nabla\times \delta H)\times H_2\big\|_{L^1_T(\tilde{B}_{2,1}^{\f12,\f32})}\\
&\qquad+\Big(1+\|c_1\|_{L^\infty_T(\dot{B}_{2,1}^{\f32}))}\Big)\big\|(\nabla\times H_1)\times \delta H\big\|_{L^1_T(\tilde{B}_{2,1}^{\f12,\f32})}\\
&\quad\lesssim \Big(1+\|c_1\|_{L^\infty_T(\dot{B}_{2,1}^{\f32}))}+\|c_2\|_{L^\infty_T(\dot{B}_{2,1}^{\f32}))}\Big)\|\delta c\|_{L^\infty_T(\tilde{B}_{2,1}^{-\f12,\f12})}\big\| H_2\big\|_{L^\infty_T(\tilde{B}_{2,1}^{\f12,\f32})}\big\|H_2\big\|_{L^1_T(\tilde{B}_{2,1}^{\f52,\f72})}\\
&\qquad+\Big(1+\|c_1\|_{L^\infty_T(\dot{B}_{2,1}^{\f32}))}\Big)\big\|\delta H\big\|_{L^2_T(\tilde{B}_{2,1}^{\f12,\f32})}\big\|H_2\big\|_{L^2_T(\tilde{B}_{2,1}^{\f32,\f52})}\\
&\qquad+\Big(1+\|c_1\|_{L^\infty_T(\dot{B}_{2,1}^{\f32}))}\Big)\big\| H_1\big\|_{L^2_T(\tilde{B}_{2,1}^{\f32,\f52})}\big\|\delta H\big\|_{L^2_T(\tilde{B}_{2,1}^{\f12,\f32})}.
 \end{split}
\end{equation*}
Substituting those estimates back into \eqref{52}, we eventually get
$$\|(\delta c, \delta u, \delta H)\|_{\mathcal{E}_{T}^{\f{1}{2}}}\le Z(T)\|(\delta c, \delta u, \delta H)\|_{\mathcal{E}_{T}^{\f{1}{2}}}$$
with
\begin{equation*}
 \begin{split}
Z(T)&=e^{\int_0^t\|u_2(\tau)\|_{\dot{B}^{\f{5}{2}}_{2,1}}d\tau}\Bigg[\|c_1\|_{L^\infty_T(\tilde{B}_{2,1}^{\f12,\f32})}+\| u_2\|_{L^1_T(\dot{B}_{2,1}^{\f52})}+\|u_2\|_{L^2_T(\dot{B}_{2,1}^{\f32})}+\| u_1\|_{L^2_T(\dot{B}_{2,1}^{\f32})}\\
&\quad+\Big(1+\|c_1\|_{L^\infty_T(\dot{B}_{2,1}^{\f32})}+\|c_2\|_{L^\infty_T(\dot{B}_{2,1}^{\f32})}\Big)\|u_2\|_{L^1_T(\dot{B}_{2,1}^{\f52})}
+T\Big(\|c_1\|_{L^\infty_T(\dot{B}_{2,1}^{\f32})}+\|c_2\|_{L^\infty_T(\dot{B}_{2,1}^{\f32})}\Big)\\
&\quad+\Big(1+\|c_1\|_{L^\infty_T(\dot{B}_{2,1}^{\f32})}+\|c_2\|_{L^\infty_T(\dot{B}_{2,1}^{\f32})}\Big)\| H_2\|_{L^\infty_T(\dot{B}_{2,1}^{\f12,\f32})}\| H_{2}\|_{L^1_T(\tilde{B}_{2,1}^{\f52,\f72})}\\
&\quad+\Big(1+\|c_1\|_{L^\infty_T(\dot{B}_{2,1}^{\f32})}\Big)\| H_{2}\|_{L^1_T(\tilde{B}_{2,1}^{\f52,\f72})}+\Big(1+\|c_1\|_{L^\infty_T(\dot{B}_{2,1}^{\f32})}\Big)\|H_2\|_{L^2_T(\tilde{B}_{2,1}^{\f32,\f52})}\\
&\quad+\| H_1\|_{L^1_T(\tilde{B}_{2,1}^{\f52,\f72})}+\| H_2\|_{L^2_T(\tilde{B}_{2,1}^{\f32,\f52})}+\| u_1\|_{L^1_T(\dot{B}_{2,1}^{\f52})}\\
&\quad+\Big(1+\|c_1\|_{L^\infty_T(\dot{B}_{2,1}^{\f32}))}+\|c_2\|_{L^\infty_T(\dot{B}_{2,1}^{\f32}))}\Big)\big\| H_2\big\|_{L^\infty_T(\tilde{B}_{2,1}^{\f12,\f32})}\big\|H_2\big\|_{L^1_T(\tilde{B}_{2,1}^{\f52,\f72})}\\
&\quad+\Big(1+\|c_1\|_{L^\infty_T(\dot{B}_{2,1}^{\f32}))}\Big)\big\|H_2\big\|_{L^2_T(\tilde{B}_{2,1}^{\f32,\f52})}
+\Big(1+\|c_1\|_{L^\infty_T(\dot{B}_{2,1}^{\f32}))}\Big)\big\| H_1\big\|_{L^2_T(\tilde{B}_{2,1}^{\f32,\f52})}\Bigg].
 \end{split}
\end{equation*}
We notice that $\limsup_{T\rightarrow 0^+}Z(T)\le
C\|c_1\|_{L^\infty(\tilde{B}^{\f12,\f32}_{2,1})}$. This is because all
other terms involve an integral in time in $L^1$ or $L^2$ sense so that as $T$ goes to zero, all those integrals will converge to
zero. Thus, if $\eta>0$ is sufficiently small, we get
$$\|(\delta c, \delta u, \delta H)\|_{\mathcal{E}_{T}^{\f{1}{2}}}=0,$$ for certain $T>0$
small enough. Thus, we have shown  uniqueness on a small time interval $[0,T]$ such that $(c_1,u_1, H_1)=(c_2,u_2, H_2)$.
This completes the proof of the uniqueness in Theorem \ref{th:main1}.

\section{Time decay estimates}
\ \ At last, we exhibit the time decay estimates of the strong
solutions to the system \eqref{equ:CMHD} for  initial data  close to a stable
equilibrium state in critical regularity framework.
 We divide it into several steps.
\subsubsection*{Step 1: Bounds for the low frequencies}
We study the following  system
\begin{align}\label{eq:low}
\left\{
\begin{aligned}
&\partial_{t}c+\dv u=f_{1},\\
&\partial_{t}u-\mathcal{A}u+\nabla c=f_{2},\\
&\partial_{t}H-\Delta H=f_{3},\\
&\dv H=0,\\
&(c,u,H)\mid_{t=0}=(c_{0},u_{0},H_{0}),
\end{aligned} \right.
\end{align}
where $f_1\eqdefa f-u\cdot\nabla c,$ $f_2\eqdefa g-u\cdot\nabla u,$ $f_3\eqdefa h-u\cdot\nabla H.$
\medbreak
Denoting by $A(D)$ the semi-group associated to \eqref{eq:low}, we have for all $q\in\mathbb{Z},$
\begin{equation}
\label{low.1}
\left(\begin{array}{c}\ddq c(t)\\\ddq u(t)\\\ddq H(t)\end{array}\right)
=e^{tA(D)}\left(\begin{array}{c}\ddq c_0\\\ddq u_0\\\ddq H_0\end{array}\right)
+\int_0^te^{(t-\tau)A(D)}\left(\begin{array}{c}\ddq f_1(\tau)\\\ddq f_2(\tau)\\\ddq f_3(\tau)\end{array}\right)d\tau.
\end{equation}
From an explicit computation of  the action of $e^{tA(D)}$ in Fourier variables (see e.g. \cite{SK}), we
discover that there exist positive constants $c_{0}$ and $C$ depending
only on $q_0$ and such that$$|\cF(e^{tA(D)} U)(\xi)|\leq Ce^{-c_0t|\xi|^2}|\cF U(\xi)|\quad\hbox{for all }\  |\xi|\leq2^{q_0}.$$
Therefore,  using Parseval's equality and the definition of $\ddq$, we get for all $q\leq q_0,$
\begin{align*}\|e^{tA(D)}\ddq U\|_{L^2}&\lesssim  e^{-\frac{c_{0}}{4}2^{2q}t}\|\ddq U\|_{L^2}.
\end{align*}
Hence, multiplying by $t^{\frac{3}{4}+\frac{s}{2}}2^{qs}$ and summing up on $q\leq q_{0}$,
\begin{equation}
\begin{split}
\label{low.2}
t^{\frac{3}{4}+\frac{s}{2}}\sum_{q\leq q_0}2^{qs}\|e^{tA(D)}\ddq U\|_{L^2}
&\lesssim
\sum_{q\leq q_0}2^{qs}e^{-\frac{c_{0}}{4}2^{2q}t}\|\ddq U\|_{L^2}t^{\frac{3}{4}+\frac{s}{2}}\\
&\lesssim
\sum_{q\leq q_0}2^{q(s+\frac{3}{2})}e^{-\frac{c_{0}}{4}2^{2q}t}\|\ddq U\|_{L^2}2^{q(-\frac{3}{2})}t^{\frac{3}{4}+\frac{s}{2}}\\
&\lesssim
\|U\|_{\dot B^{-\frac{3}{2}}_{2,\infty}}^\ell\sum_{q\leq q_0}2^{q(s+\frac{3}{2})}e^{-\frac{c_{0}}{4}2^{2q}t}t^{\frac{3}{4}+\frac{s}{2}}\\
&\lesssim
\|U\|_{\dot B^{-\frac{3}{2}}_{2,\infty}}^\ell\sum_{q\leq q_0}2^{q(s+\frac{3}{2})}e^{-\frac{c_{0}}{4}2^{2q}t}
t^{\frac{1}{2}(s+\frac{3}{2})}.
\end{split}
\end{equation}
As for any $\sigma>0$ there  exists a constant $C_\sigma$ so that
\begin{equation}\label{low.3}
\sup_{t\geq0}\sum_{q\in\mathbb{Z}}t^{\frac\sigma2}2^{q\sigma}e^{-\frac{c_0}{4}2^{2q}t}\leq C_\sigma.
\end{equation}
We get from \eqref{low.2} and \eqref{low.3} that for $s>-3/2,$
$$
\sup_{t\geq0}\, t^{\frac 34+\frac s2}\|e^{tA(D)}U\|_{\dot B^s_{2,1}}^\ell
\lesssim\|U\|_{\dot B^{-\frac{3}{2}}_{2,\infty}}^\ell.
$$
It is also obvious that  for $s>-3/2,$
$$
\|e^{tA(D)}U\|_{\dot B^s_{2,1}}^\ell
\lesssim \|U\|_{\dot B^{-\frac{3}{2}}_{2,\infty}}^\ell\sum_{q\leq q_0}2^{q(s+\frac{3}{2})}\lesssim\|U\|_{\dot B^{-\frac{3}{2}}_{2,\infty}}^\ell .
$$
So, setting $\langle t\rangle\eqdefa\sqrt{1+t^{2}}$, we arrive at
\begin{equation}
\label{U}
\sup_{t\geq0}\, \langle t\rangle^{\frac 34+\frac s2}\|e^{tA(D)}U\|_{\dot B^s_{2,1}}^\ell
\lesssim\|U\|_{\dot B^{-\frac{3}{2}}_{2,\infty}}^\ell,
\end{equation}
and thus, taking advantage of Duhamel's formula,
\begin{equation}\label{low.4}
\Big\|\int_0^te^{(t-\tau)A(D)}(f_{1},f_{2},f_{3})(\tau)d\tau\Big\|_{\dot B^{s}_{2,1}}^\ell
\lesssim\int_0^t\langle t-\tau\rangle^{-(\frac{3}{4}+\frac{s}{2})} \big\|(f_{1},f_{2},f_{3})(\tau)\big\|_{\dot B^{-\frac{3}{2}}_{2,\infty}}^\ell d\tau.
\end{equation}
We claim that for all $s\in(-3/2,2]$,  then we have for all $t\geq0$,
\begin{equation}\label{f123low1}
\int_0^t\langle t-\tau\rangle^{-(\frac{3}{4}+\frac{s}{2})} \big\|(f_{1},f_{2},f_{3})(\tau)\big\|_{\dot B^{-\frac{3}{2}}_{2,\infty}}^\ell d\tau
\lesssim\langle t\rangle^{-(\frac 34+\frac s2)}
\big(X^2(t)+D^2(t)+D^3(t)\big).
\end{equation}
 Owing to the embedding $L^{1}\hookrightarrow\dot B^{-\frac{3}{2}}_{2,\infty}$,  it suffices to prove \eqref{f123low1}  with $\|(f_{1},f_{2},f_{3})(\tau)\|_{L^1}^{\ell}$ instead of $\|(f_{1},f_{2},f_{3})(\tau)\|_{\dot  B^{-\frac{3}{2}}_{2,\infty}}^{\ell}$.\medbreak
To bound the term with $f_1$, we use the following decomposition:
$$
f_1=u\cdot\nabla c+ c\, \div u^\ell + c\,\div u^h.
$$
Now, from H\"{o}lder's inequality,  the embedding $\dot B^{0}_{2,1}\hookrightarrow L^{2}$, the definitions of $D(t), \alpha$ and Lemma \ref{lemma2.13}, one may write for all $s\in(-\frac 32,2]$,
\begin{equation}
\begin{split}\label{udivc}
&\int_0^t\langle t-\tau\rangle^{-(\frac 34+\frac s2)}\|(u\cdot\nabla c)(\tau)\|_{L^1}\,d\tau\\
&\quad\lesssim\int_0^t\langle t-\tau\rangle^{-(\frac 34+\frac s2)}\|u\|_{L^2}\|\nabla c\|_{L^2}\,d\tau\\
&\quad\lesssim\int_0^t\langle t-\tau\rangle^{-(\frac 34+\frac s2)}\|u\|_{\dot  B^{0}_{2,1}}\|\nabla c\|_{\dot  B^{0}_{2,1}}\,d\tau\\
&\quad\lesssim\int_0^t\langle t-\tau\rangle^{-(\frac 34+\frac s2)}\big(\|u\|_{\dot  B^{0}_{2,1}}^{\ell}+\|u\|_{\dot  B^{0}_{2,1}}^h\big)\big(\|\nabla c\|_{\dot  B^{0}_{2,1}}^{\ell}+\|\nabla c\|_{\dot  B^{0}_{2,1}}^{h}\big)\,d\tau\\
&\quad\lesssim\big(\sup_{0\leq\tau\leq t}\langle \tau\rangle^{\frac 34}\|u(\tau)\|_{\dot  B^{0}_{2,1}}^{\ell}\big)
\big(\sup_{0\leq\tau\leq t}\langle \tau\rangle^{\frac 54}\|\nabla c(\tau)\|_{\dot  B^{0}_{2,1}}^{\ell}\big)
\int_0^t\langle t-\tau\rangle^{-(\frac 34+\frac s2)}\langle \tau\rangle^{-\frac 34}\langle\tau\rangle^{-\frac 54}\,d\tau\\
&\qquad+\big(\sup_{0\leq\tau\leq t}\langle \tau\rangle^{\frac 34}\|u(\tau)\|_{\dot  B^{0}_{2,1}}^{\ell}\big)
\big(\sup_{0\leq\tau\leq t}\langle \tau\rangle^\alpha\|\nabla c(\tau)\|_{\dot  B^{\frac{1}{2}}_{2,1}}^h\big)
\int_0^t\langle t-\tau\rangle^{-(\frac 34+\frac s2)}\langle \tau\rangle^{-\frac 34}\langle\tau\rangle^{-\alpha}\,d\tau\\
&\qquad+\big(\sup_{0\leq\tau\leq t}\langle \tau\rangle^\alpha\|u(\tau)\|_{\dot  B^{\frac{1}{2}}_{2,1}}^h\big)
\big(\sup_{0\leq\tau\leq t}\langle \tau\rangle^{\frac 54}\|\nabla c(\tau)\|_{\dot  B^{0}_{2,1}}^{\ell}\big)
\int_0^t\langle t-\tau\rangle^{-(\frac 34+\frac s2)}\langle \tau\rangle^{-\alpha}\langle\tau\rangle^{-\frac 54}\,d\tau\\
&\qquad+\big(\sup_{0\leq\tau\leq t}\langle \tau\rangle^\alpha\|u(\tau)\|_{\dot  B^{\frac{1}{2}}_{2,1}}^h\big)
\big(\sup_{0\leq\tau\leq t}\langle \tau\rangle^{\alpha}\|\nabla c(\tau)\|_{\dot  B^{\frac{1}{2}}_{2,1}}^h\big)
\int_0^t\langle t-\tau\rangle^{-(\frac 34+\frac s2)}\langle \tau\rangle^{-\alpha}\langle\tau\rangle^{-\alpha}\,d\tau\\
&\quad\lesssim D^{2}(t)\int_0^t\langle t-\tau\rangle^{-(\frac 34+\frac s2)}\langle \tau\rangle^{-\min({2,\alpha+\frac{3}{4},2\alpha})}\,d\tau\\
&\quad\lesssim \langle t\rangle^{-(\frac 34+\frac s2)} D^{2}(t).
\end{split}
\end{equation}
The term  $c\,\div u^\ell$ may be treated along the same lines, we have
\begin{equation}
\begin{split}
\label{cdivul}
&\int_0^t\langle t-\tau\rangle^{-(\frac 34+\frac s2)}\|c\div u^{\ell} \|_{L^1}\,d\tau\\
&\quad\lesssim\int_0^t\langle t-\tau\rangle^{-(\frac 34+\frac s2)}\big(\|c\|_{\dot  B^{0}_{2,1}}^{\ell}+\|c\|_{\dot  B^{0}_{2,1}}^h\big)\|\nabla u\|_{\dot  B^{0}_{2,1}}^{\ell}\,d\tau\\
&\quad\lesssim\big(\sup_{0\leq\tau\leq t}\langle \tau\rangle^{\frac 34}\|c(\tau)\|_{\dot  B^{0}_{2,1}}^{\ell}\big)
\big(\sup_{0\leq\tau\leq t}\langle \tau\rangle^\frac{5}{4}\| u(\tau)\|_{\dot  B^{1}_{2,1}}^\ell\big)
\int_0^t\langle t-\tau\rangle^{-(\frac 34+\frac s2)}\langle \tau\rangle^{-\frac 34}\langle\tau\rangle^{-\frac{5}{4}}\,d\tau\\
&\qquad+\big(\sup_{0\leq\tau\leq t}\langle \tau\rangle^\alpha\|c(\tau)\|_{\dot  B^{\frac{3}{2}}_{2,1}}^h\big)
\big(\sup_{0\leq\tau\leq t}\langle \tau\rangle^{\frac 54}\| u(\tau)\|_{\dot  B^{1}_{2,1}}^{\ell}\big)
\int_0^t\langle t-\tau\rangle^{-(\frac 34+\frac s2)}\langle \tau\rangle^{-\alpha}\langle\tau\rangle^{-\frac 54}\,d\tau\\
&\quad\lesssim D^{2}(t)\int_0^t\langle t-\tau\rangle^{-(\frac 34+\frac s2)}\langle \tau\rangle^{-\min({2,\alpha+\frac{5}{4}})}\,d\tau\\
&\quad\lesssim \langle t\rangle^{-(\frac 34+\frac s2)} D^{2}(t).
\end{split}
\end{equation}
Regarding the term with $c\,\div u^h,$ we use that if $t\geq2,$
\begin{equation*}
\begin{split}
&\int_0^t\langle t-\tau\rangle^{-(\frac 34+\frac s2)}\|(c\div u^h)(\tau)\|_{L^1}\,d\tau\\
&\quad\lesssim\int_0^t\langle t-\tau\rangle^{-(\frac 34+\frac s2)} \|c(\tau)\|_{L^2}\|\div u(\tau)\|^h_{L^2}\,d\tau\\
&\quad\lesssim\int_0^t\langle t-\tau\rangle^{-(\frac 34+\frac s2)} \|c(\tau)\|_{\dot  B^{0}_{2,1}}\|\div u(\tau)\|^h_{\dot  B^{0}_{2,1}}\,d\tau\\
&\quad\lesssim\int_0^1\langle t-\tau\rangle^{-(\frac 34+\frac s2)} \|c(\tau)\|_{\dot  B^{0}_{2,1}}\|\div u(\tau)\|^h_{\dot  B^{0}_{2,1}}\,d\tau\\
&\qquad+\int_1^t\langle t-\tau\rangle^{-(\frac 34+\frac s2)} \|c(\tau)\|_{\dot  B^{0}_{2,1}}\|\div u(\tau)\|^h_{\dot  B^{0}_{2,1}}\,d\tau\\
&\quad\eqdefa I_{1}+I_{2}.
\end{split}
\end{equation*}
Remembering the definitions of $X(t)$ and $D(t)$, we  obtain
\begin{align*}
I_{1}&=\int_0^1\langle t-\tau\rangle^{-(\frac 34+\frac s2)} \|c(\tau)\|_{\dot  B^{0}_{2,1}}\|\div u(\tau)\|^h_{\dot  B^{0}_{2,1}}\,d\tau\\
&\lesssim\langle t\rangle^{-(\frac 34+\frac s2)} \sup _{0\leq\tau\leq1}\|c(\tau)\|_{\dot  B^{0}_{2,1}}\int_0^1\|\div u(\tau)\|^h_{\dot  B^{0}_{2,1}}\,d\tau\\
&\lesssim\langle t\rangle^{-(\frac 34+\frac s2)} \sup _{0\leq\tau\leq1}\|c(\tau)\|_{\dot  B^{0}_{2,1}}\int_0^1\| u(\tau)\|^h_{\dot  B^{\frac{5}{2}}_{2,1}}\,d\tau\\
&\lesssim\langle t\rangle^{-(\frac 34+\frac s2)} D(1)X(1)
\end{align*}
and,  using the fact that $\langle \tau\rangle\approx\tau$ when $\tau\geq1$,
\begin{align*}
I_{2} &=\int_1^t\langle t-\tau\rangle^{-(\frac 34+\frac s2)} \|c(\tau)\|_{\dot  B^{0}_{2,1}}\|\div u(\tau)\|^h_{\dot  B^{0}_{2,1}}\,d\tau\\
&\lesssim\int_1^t\langle t-\tau\rangle^{-(\frac 34+\frac s2)}
\big(\|c(\tau)\|_{\dot  B^{0}_{2,1}}^{\ell}+\|c(\tau)\|_{\dot  B^{0}_{2,1}}^{h}\big)\|\div u(\tau)\|^h_{\dot  B^{0}_{2,1}}\,d\tau\\
&\lesssim\big(\sup_{1\leq\tau\leq t}\langle \tau\rangle^{\frac 34}\|c(\tau)\|_{\dot  B^{0}_{2,1}}^{\ell}\big)
\big(\sup_{1\leq\tau\leq t}\| \tau\nabla u(\tau)\|^{h}_{\dot  B^{\frac{3}{2}}_{2,1}}\big)
\int_1^t\langle t-\tau\rangle^{-(\frac 34+\frac s2)}\langle \tau\rangle^{-\frac 34}\langle\tau\rangle^{-1}\,d\tau\\
&\quad+\big(\sup_{1\leq\tau\leq t}\langle \tau\rangle^{\alpha}\|c(\tau)\|_{\dot  B^{\frac{3}{2}}_{2,1}}^h\big)
\big(\sup_{1\leq\tau\leq t}\| \tau\nabla u(\tau)\|^{h}_{\dot  B^{\frac{3}{2}}_{2,1}}\big)
\int_1^t\langle t-\tau\rangle^{-(\frac 34+\frac s2)}\langle \tau\rangle^{-\alpha}\langle\tau\rangle^{-1}\,d\tau\\
&\lesssim D^{2}(t)\int_1^t\langle t-\tau\rangle^{-(\frac 34+\frac s2)}\langle \tau\rangle^{-\min({\alpha+1,\frac{7}{4}})}\,d\tau\\
&\lesssim \langle t\rangle^{-(\frac 34+\frac s2)} D^{2}(t).
\end{align*}
Therefore, for $t\geq2,$ we arrive at
\begin{align}\label{cdivuh}
&\int_0^t\langle t-\tau\rangle^{-(\frac 34+\frac s2)}\|c\div u^h(\tau)\|_{L^1}\,d\tau\nonumber\\
&\quad\lesssim\langle t\rangle^{-(\frac 34+\frac s2)} \big(D^{2}(t)+X^{2}(t)\big).
\end{align}
The case $t\leq2$ is obvious as $\langle t\rangle\approx1$ and
$\langle t-\tau\rangle\approx1$ for $0\leq\tau\leq t\leq 2$, and
\begin{equation}
\begin{split}
\label{cdivuht}
&\int_0^t\|c\,\div u^h\|_{L^1}\,d\tau\\
&\quad\lesssim \|c\|_{L^\infty_t(L^2)}\|\div u\|^h_{L_t^1(L^2)}\\
&\quad\lesssim \|c\|_{L^\infty_t(\dot  B^{0}_{2,1})}\|\div u\|^h_{L_t^1(\dot  B^{0}_{2,1})}\\
&\quad\lesssim \|c\|_{L^\infty_t(\dot  B^{0}_{2,1})}\| u\|^h_{L_t^1(\dot  B^{\frac{5}{2}}_{2,1})}\\
&\quad\lesssim X(t)D(t).
\end{split}
\end{equation}
From \eqref{udivc}-\eqref{cdivuht}, we get
\begin{equation*}\label{f123low}
\int_0^t\langle t-\tau\rangle^{-(\frac{3}{4}+\frac{s}{2})} \big\|f_{1}(\tau)\big\|_{\dot B^{-\frac{3}{2}}_{2,\infty}}^\ell d\tau
\lesssim\langle t\rangle^{-(\frac 34+\frac s2)}
\big(X^2(t)+D^2(t)\big).
\end{equation*}
Next, in order to bound the term of \eqref{f123low} corresponding to $f_2$, we use the following decomposition
\begin{align*}
f_2&=g-u\cdot\nabla u\\
&=-u\cdot\nabla u^{\ell}-u\cdot\nabla u^{h}-L_{1}(c)\cA u+ L_{2}(c)\nabla c\\
&\quad- L_{3}(c)\big(\frac{1}{2}\nabla|H|^{2}-H\cdot\nabla H\big).
\end{align*}
Similar to  \eqref{cdivul}-\eqref{cdivuht}, we have
$$\int_0^t\langle t-\tau\rangle^{-(\frac 34+\frac s2)}\|u\cdot\nabla u^{\ell}\|_{L^1}\,d\tau\lesssim \langle t\rangle^{-(\frac 34+\frac s2)} D^{2}(t)$$
and
\begin{align*}
\int_0^t\langle t-\tau\rangle^{-(\frac 34+\frac s2)}\|u\cdot\nabla u^h(\tau)\|_{L^1}\,d\tau
\lesssim\langle t\rangle^{-(\frac 34+\frac s2)} \big(D^{2}(t)+X^{2}(t)\big).
\end{align*}
For $L_{1}(c)\cA u,$ we write that
$$
L_{1}(c)\cA u=L_{1}(c)\cA u^\ell+ L_{1}(c)\cA u^h,
$$
where $L_{1}$ stands for some smooth function vanishing at $0$.
Now, we have
 \begin{align*}
 &\int_0^t\langle t-\tau\rangle^{-(\frac 34+\frac s2)}\|L_{1}(c)\cA u^\ell\|_{L^1}\,d\tau\\
 &\quad\lesssim\int_0^t\langle t-\tau\rangle^{-(\frac 34+\frac s2)}\big(\|c\|_{\dot  B^{0}_{2,1}}^{\ell}+\|c\|_{\dot  B^{0}_{2,1}}^h\big)\|\cA u^\ell\|_{\dot  B^{0}_{2,1}}\,d\tau\\
 &\quad\lesssim \big(\sup_{\tau\in[0,t]} \langle\tau\rangle^{\frac 34}\|c(\tau)\|_{\dot  B^{0}_{2,1}}^{l}\big)
 \big(\sup_{\tau\in[0,t]} \langle\tau\rangle^{\frac 74}\|\cA u^\ell\|_{\dot  B^{0}_{2,1}}\big)
 \int_0^t \langle t-\tau\rangle^{-(\frac 34+\frac s2)}\langle\tau\rangle^{-\frac 52}\,d\tau\\
 &\qquad+ \big(\sup_{1\leq\tau\leq t}\langle \tau\rangle^{\alpha}\|c(\tau)\|_{\dot  B^{\frac{3}{2}}_{2,1}}^h\big)
  \big(\sup_{\tau\in[0,t]} \langle\tau\rangle^{\frac 74}\|\cA u^\ell\|_{\dot  B^{0}_{2,1}}\big)
 \int_0^t \langle t-\tau\rangle^{-(\frac 34+\frac s2)}\langle\tau\rangle^{-(\frac 74+\alpha)}\,d\tau\\
 &\quad\lesssim D^{2}(t)\int_0^t\langle t-\tau\rangle^{-(\frac 34+\frac s2)}\langle \tau\rangle^{-\min(\frac 52,\frac 74+\alpha)}\,d\tau\\
 &\quad\lesssim \langle t\rangle^{-(\frac 34+\frac s2)} D^{2}(t).
 \end{align*}
To handle the term $L_{1}(c)\cA u^h$, we consider the cases $t\geq2$ and $t\leq2$ separately. If $t\geq2$ then we write
 \begin{align*}
 &\int_0^t\langle t-\tau\rangle^{-(\frac 34+\frac s2)}\|L_{1}(c)\cA u^h\|_{L^1}\,d\tau
 \lesssim\int_0^t\langle t-\tau\rangle^{-(\frac 34+\frac s2)}\|c\|_{\dot  B^{0}_{2,1}}\|\cA u\|^h_{\dot  B^{0}_{2,1}}\,d\tau\\
 &\quad\lesssim\int_0^1\langle t-\tau\rangle^{-(\frac 34+\frac s2)}\|c\|_{\dot  B^{0}_{2,1}}\|\cA u\|^h_{\dot  B^{0}_{2,1}}\,d\tau
 +\int_1^t\langle t-\tau\rangle^{-(\frac 34+\frac s2)}\|c\|_{\dot  B^{0}_{2,1}}\|\cA u\|^h_{\dot  B^{0}_{2,1}}\,d\tau\\
 &\quad\eqdefa K_{1}+K_{2}.
 \end{align*}
From the definitions of $X(t)$ and $D(t)$, we  obtain
\begin{align*}
 K_{1}&=\int_0^1\langle t-\tau\rangle^{-(\frac 34+\frac s2)}\|c\|_{\dot  B^{0}_{2,1}}\|\cA u\|^h_{\dot  B^{0}_{2,1}}\,d\tau\\
 &\lesssim\langle t\rangle^{-(\frac 34+\frac s2)} \big(\sup_{\tau\in[0,1]}\|c(\tau)\|_{\dot  B^{0}_{2,1}}\big)
\int_0^1 \|u\|^h_{\dot  B^{\frac{5}{2}}_{2,1}}\,d\tau\\
&\lesssim\langle t\rangle^{-(\frac 34+\frac s2)} D(1)X(1),
\end{align*}
and, using the fact that $\langle \tau\rangle\approx\tau$ when $\tau\geq1$,
\begin{align*}
 K_{2}&=\int_1^t\langle t-\tau\rangle^{-(\frac 34+\frac s2)}\|c\|_{\dot  B^{0}_{2,1}}\|\cA u\|^h_{\dot  B^{0}_{2,1}}\,d\tau\\
&\lesssim\int_1^t\langle t-\tau\rangle^{-(\frac 34+\frac s2)}\big(\|c\|_{\dot  B^{0}_{2,1}}^{\ell}+\|c\|_{\dot  B^{0}_{2,1}}^{h}\big)\|\cA u\|^h_{\dot  B^{0}_{2,1}}\,d\tau\\
&\lesssim\int_1^t\langle t-\tau\rangle^{-(\frac 34+\frac s2)}\|c\|_{\dot  B^{0}_{2,1}}^{\ell}\|\cA u\|^h_{\dot  B^{0}_{2,1}}\,d\tau
+\int_1^t\langle t-\tau\rangle^{-(\frac 34+\frac s2)}\|c\|_{\dot  B^{0}_{2,1}}^{h})\|\cA u\|^h_{\dot  B^{0}_{2,1}}\,d\tau\\
&\lesssim\big(\sup_{0\leq\tau\leq t}\langle \tau\rangle^{\frac 34}\|c(\tau)\|_{\dot  B^{0}_{2,1}}^{\ell}\big)
\big(\sup_{0\leq\tau\leq t} \|\tau\nabla u(\tau)\|_{\dot  B^{\frac{3}{2}}_{2,1}}^{h}\big)
\int_1^t\langle t-\tau\rangle^{-(\frac 34+\frac s2)}\langle \tau\rangle^{-\frac 74}\,d\tau\\
&+\big(\sup_{0\leq\tau\leq t}\langle \tau\rangle^{\alpha}\|\nabla c(\tau)\|_{\dot  B^{\frac{1}{2}}_{2,1}}^h\big)
\big(\sup_{0\leq\tau\leq t} \|\tau\nabla u(\tau)\|_{\dot  B^{\frac{3}{2}}_{2,1}}^h\big)
\int_1^t\langle t-\tau\rangle^{-(\frac 34+\frac s2)}\langle \tau\rangle^{-(\alpha+1)}d\tau\\
&\lesssim D^{2}(t)
\int_1^t\langle t-\tau\rangle^{-(\frac 34+\frac s2)}\langle \tau\rangle^{-\min(\alpha+1,\frac 74)}d\tau\\
&\lesssim \langle t\rangle^{-(\frac 34+\frac s2)} D^{2}(t).
\end{align*}
Thus, for $t\geq2,$ we arrive at
$$\int_0^t\langle t-\tau\rangle^{-(\frac 34+\frac s2)}\|L_{1}(c)\cA u^h\|_{L^1}\,d\tau
\lesssim\langle t\rangle^{-(\frac 34+\frac s2)}(X^{2}(t)+D^{2}(t)).$$
The case $t\leq2$ is obvious as $\langle t\rangle\approx1$ and $\langle t-\tau\rangle\approx1$ for $0\leq\tau\leq t\leq2,$
\begin{align*}
&\int_0^t\|L_{1}(c)\cA u^h\|_{L^1}\,d\tau\\
&\quad\lesssim\int_0^t\|c\|_{\dot  B^{0}_{2,1}}\|\cA u^h\|_{\dot  B^{0}_{2,1}}\,d\tau\\
&\quad\lesssim\big(\sup_{\tau\in[0,1]}\|c(\tau)\|_{\dot  B^{0}_{2,1}}\big)
\int_0^1 \|u^h\|_{\dot  B^{\frac{5}{2}}_{2,1}}\,d\tau\\
&\quad\lesssim D(t)X(t).
\end{align*}
Similar to  \eqref{udivc},  we have
$$\int_0^t\langle t-\tau\rangle^{-(\frac 34+\frac s2)}\|L_{2}(c)\nabla c\|_{L^1}\,d\tau\lesssim \langle t\rangle^{-(\frac 34+\frac s2)} D^{2}(t).$$
To Bound the term $L_{3}(c)\big(\frac{1}{2}\nabla|H|^{2}-H\cdot\nabla H\big)$, it suffices to consider  $L_{3}(c)H\cdot\nabla H$ as follows
\begin{align*}
&\int_0^t\langle t-\tau\rangle^{-(\frac 34+\frac s2)}\Big\|L_{3}(c)\big(\frac{1}{2}\nabla|H|^{2}-H\cdot\nabla H\big)\Big\|_{L^1}\,d\tau\\
&\quad\lesssim\int_0^t\langle t-\tau\rangle^{-(\frac 34+\frac s2)}\|L_{3}(c)\|_{L^2}\|H\cdot\nabla H\|_{L^2}\,d\tau\\
&\quad\lesssim\int_0^t\langle t-\tau\rangle^{-(\frac 34+\frac s2)}\|H\|_{\dot  B^{\frac{3}{2}}_{2,1}}\|\nabla H\|_{\dot  B^{0}_{2,1}}d\tau\\
&\qquad+\int_0^t\langle t-\tau\rangle^{-(\frac 34+\frac s2)}\|c\|_{\dot  B^{0}_{2,1}}\|H\|_{\dot  B^{\frac{3}{2}}_{2,1}}\|\nabla H\|_{\dot  B^{0}_{2,1}}\,d\tau\\
 &\quad\eqdefa L_{1}+L_{2}.
\end{align*}
We estimate the two terms  $L_{1}$ and $L_{2}$ in the following,
\begin{align*}
L_{1}&\lesssim\int_0^t\langle t-\tau\rangle^{-(\frac 34+\frac s2)}(\|H\|_{\dot  B^{\frac{3}{2}}_{2,1}}^{\ell}+\|H\|_{\dot  B^{\frac{3}{2}}_{2,1}}^{h})
(\|\nabla H\|_{\dot  B^{0}_{2,1}}^{\ell}+\|\nabla H\|_{\dot  B^{0}_{2,1}}^{h}d\tau)\\
&\lesssim\big(\sup_{0\leq\tau\leq t}\langle \tau\rangle^{\frac 32}\|H(\tau)\|_{\dot  B^{\frac{3}{2}}_{2,1}}^{\ell}\big)
\big(\sup_{0\leq\tau\leq t}\langle \tau\rangle^\frac{5}{4}\| H(\tau)\|_{\dot  B^{1}_{2,1}}^\ell\big)
\int_0^t\langle t-\tau\rangle^{-(\frac 34+\frac s2)}\langle \tau\rangle^{-\frac 32}\langle\tau\rangle^{-\frac{5}{4}}\,d\tau\\
&\quad+\big(\sup_{0\leq\tau\leq t}\langle \tau\rangle^{\frac 32}\|H(\tau)\|_{\dot  B^{\frac{3}{2}}_{2,1}}^{\ell}\big)
\big(\sup_{0\leq\tau\leq t}\langle \tau\rangle^\alpha\|\nabla H(\tau)\|_{\dot  B^{\frac{1}{2}}_{2,1}}^h\big)
\int_0^t\langle t-\tau\rangle^{-(\frac 34+\frac s2)}\langle\tau\rangle^{-(\alpha
+\frac 32)}\,d\tau\\
&\quad+\big(\sup_{0\leq\tau\leq t}\langle \tau\rangle^\alpha\| H(\tau)\|_{\dot  B^{\frac{3}{2}}_{2,1}}^h\big)
\big(\sup_{0\leq\tau\leq t}\langle \tau\rangle^\frac{5}{4}\| H(\tau)\|_{\dot  B^{1}_{2,1}}^\ell\big)
\int_0^t\langle t-\tau\rangle^{-(\frac 34+\frac s2)}\langle \tau\rangle^{-(\alpha+\frac{5}{4})}\,d\tau\\
&\quad+\big(\sup_{0\leq\tau\leq t}\langle \tau\rangle^\alpha\| H(\tau)\|_{\dot  B^{\frac{3}{2}}_{2,1}}^h\big)
\big(\sup_{0\leq\tau\leq t}\langle \tau\rangle^\alpha\|\nabla H(\tau)\|_{\dot  B^{\frac{1}{2}}_{2,1}}^h\big)
\int_0^t\langle t-\tau\rangle^{-(\frac 34+\frac s2)}\langle \tau\rangle^{-2\alpha}\,d\tau\\
&\lesssim D^{2}(t)\int_0^t\langle t-\tau\rangle^{-(\frac 34+\frac s2)}\langle \tau\rangle^{-\min(\frac{11}{4},\alpha+\frac{5}{4},2\alpha)}\,d\tau\\
&\lesssim \langle t\rangle^{-(\frac 34+\frac s2)} D^{2}(t),
\end{align*} and
\begin{align*}
L_{2}
&\lesssim\int_0^t\langle t-\tau\rangle^{-(\frac 34+\frac s2)}\big(\|c\|_{\dot  B^{0}_{2,1}}^{\ell}+\|c\|_{\dot  B^{0}_{2,1}}^{h}\big)\|H\|_{\dot  B^{\frac{3}{2}}_{2,1}}\|\nabla H\|_{\dot  B^{0}_{2,1}}\,d\tau\\
&\lesssim\big(\sup_{0\leq\tau\leq t}\langle \tau\rangle^{\frac 34}\|c(\tau)\|_{\dot  B^{0}_{2,1}}^{\ell}\big)
\int_0^t\langle t-\tau\rangle^{-(\frac 34+\frac s2)}\langle \tau\rangle^{-\frac 34}\|H\|_{\dot  B^{\frac{3}{2}}_{2,1}}\|\nabla H\|_{\dot  B^{0}_{2,1}}\,d\tau\\
&\quad+\big(\sup_{0\leq\tau\leq t}\langle \tau\rangle^\alpha\|c(\tau)\|_{\dot  B^{\frac{3}{2}}_{2,1}}^h\big)
\int_0^t\langle t-\tau\rangle^{-(\frac 34+\frac s2)}\langle \tau\rangle^{-\alpha}\|H\|_{\dot  B^{\frac{3}{2}}_{2,1}}\|\nabla H\|_{\dot  B^{0}_{2,1}}\,d\tau\\
&\lesssim D^{3}(t)\int_0^t\langle t-\tau\rangle^{-(\frac 34+\frac s2)}\langle \tau\rangle^{-\min(\frac{7}{2},\alpha+2,2\alpha+\frac 34,3\alpha)}\,d\tau\\
&\lesssim \langle t\rangle^{-(\frac 34+\frac s2)} D^{3}(t).
\end{align*}
Thus,
\begin{equation*}\label{f123low}
\int_0^t\langle t-\tau\rangle^{-(\frac{3}{4}+\frac{s}{2})} \big\|f_{2}(\tau)\big\|_{\dot B^{-\frac{3}{2}}_{2,\infty}}^\ell d\tau
\lesssim\langle t\rangle^{-(\frac 34+\frac s2)}
\big(X^2(t)+D^2(t)+D^3(t)\big).
\end{equation*}
Finally, we bound the term  $f_3$. Similar to \eqref{udivc}-\eqref{cdivuht},
we have
\begin{align*}&\int_0^t\langle t-\tau\rangle^{-(\frac 34+\frac s2)}\big(\|u\cdot\nabla H\|_{L^1}+\|H(\div u)\|_{L^1}+\|H\cdot\nabla u\|_{L^1}\big)\,d\tau\\
&\quad\lesssim \langle t\rangle^{-(\frac 34+\frac s2)}\big(X^{2}(t)+D^{2}(t)\big).
\end{align*}
To bound the term $\nabla\times\big(L_{3}(c)(\nabla\times H)\times H\big)$,  employing Bernstein's inequality and Proposition \ref{p26}, we conclude that
\begin{align*}
&\int_0^t\langle t-\tau\rangle^{-(\frac 34+\frac s2)} \Big\|\nabla\times\big(L_{3}(c)(\nabla\times H)\times H\big)\Big\|_{\dot  B^{-\frac{3}{2}}_{2,\infty}}\,d\tau\\
&\quad\lesssim\int_0^t\langle t-\tau\rangle^{-(\frac 34+\frac s2)} \Big\|L_{3}(c)(\nabla\times H)\times H\Big\|_{\dot  B^{-\frac{1}{2}}_{2,\infty}}\,d\tau\\
&\quad\lesssim\int_0^t\langle t-\tau\rangle^{-(\frac 34+\frac s2)}
\| H\|_{\dot  B^{\frac{3}{2}}_{2,1}}
\|\nabla H\|_{\dot  B^{-\frac{1}{2}}_{2,1}}\,d\tau\\
&\qquad+\int_0^t\langle t-\tau\rangle^{-(\frac 34+\frac s2)}
\|c\|_{\dot  B^{\frac{3}{2}}_{2,1}}\| H\|_{\dot  B^{\frac{3}{2}}_{2,1}}
\|\nabla H\|_{\dot  B^{-\frac{1}{2}}_{2,1}}\,d\tau\\
&\quad\eqdefa M_{1}+M_{2}.
\end{align*}
We bound  the terms  $M_{1}$ and $M_{2}$ as follows
\begin{align*}
M_{1}
&\lesssim\int_0^t\langle t-\tau\rangle^{-(\frac 34+\frac s2)}
\big(\| H\|_{\dot  B^{\frac{3}{2}}_{2,1}}^{\ell}
+\| H\|_{\dot  B^{\frac{3}{2}}_{2,1}}^{h}\big)
\big(\|\nabla H\|_{\dot  B^{-\frac{1}{2}}_{2,1}}^{\ell}
+\|\nabla H\|_{\dot  B^{-\frac{1}{2}}_{2,1}}^{h}\big)\,d\tau\\
&\lesssim
\big(\sup_{0\leq\tau\leq t}\langle \tau\rangle^\frac{3}{2}\| H(\tau)\|_{\dot  B^{\frac{3}{2}}_{2,1}}^\ell\big)
\big(\sup_{0\leq\tau\leq t}\langle \tau\rangle\| H(\tau)\|_{\dot  B^{\frac{1}{2}}_{2,1}}^\ell\big)
\int_0^t\langle t-\tau\rangle^{-(\frac 34+\frac s2)}\langle \tau\rangle^{-\frac{5}{2}}\,d\tau\\
&\quad+\big(\sup_{0\leq\tau\leq t}\langle \tau\rangle^\frac{3}{2}\|H(\tau)\|_{\dot  B^{\frac{3}{2}}_{2,1}}^{\ell}\big)
\big(\sup_{0\leq\tau\leq t}\langle \tau\rangle^\alpha\|\nabla H(\tau)\|_{\dot  B^{\frac{1}{2}}_{2,1}}^{h}\big)
\int_0^t\langle t-\tau\rangle^{-(\frac 32+\alpha)}\langle \tau\rangle^{-(\alpha+\frac{3}{2})}\,d\tau\\
&\quad+\big(\sup_{0\leq\tau\leq t}\langle \tau\rangle^\alpha
\| H(\tau)\|_{\dot  B^{\frac{3}{2}}_{2,1}}^h\big)
\big(\sup_{0\leq\tau\leq t}\langle \tau\rangle\| H(\tau)\|_{\dot  B^{\frac{1}{2}}_{2,1}}^\ell\big)
\int_0^t\langle t-\tau\rangle^{-(\frac 34+\frac s2)}\langle \tau\rangle^{-(\alpha+1)}\,d\tau\\
&\quad+\big(\sup_{0\leq\tau\leq t}\langle \tau\rangle^\alpha\|H(\tau)\|_{\dot  B^{\frac{3}{2}}_{2,1}}^{h}\big)
\big(\sup_{0\leq\tau\leq t}\langle \tau\rangle^\alpha\|\nabla H(\tau)\|_{\dot  B^{\frac{1}{2}}_{2,1}}^{h}\big)
\int_0^t\langle t-\tau\rangle^{-(\frac 34+\frac s2)}\langle \tau\rangle^{-2\alpha}d\tau\\
&\lesssim D^{2}(t)\int_0^t\langle t-\tau\rangle^{-(\frac 34+\frac s2)}\langle \tau\rangle^{-\min(\frac{5}{2},\alpha+1,2\alpha)}\,d\tau\\
&\lesssim \langle t\rangle^{-(\frac 34+\frac s2)} D^{2}(t),
\end{align*}
\begin{align*}
M_{2}&\lesssim\int_0^t\langle t-\tau\rangle^{-(\frac 34+\frac s2)}\big(
\|c\|_{\dot  B^{\frac{3}{2}}_{2,1}}^{\ell}+\|c\|_{\dot  B^{\frac{3}{2}}_{2,1}}^{h}\big)\| H\|_{\dot  B^{\frac{3}{2}}_{2,1}}
\|\nabla H\|_{\dot  B^{-\frac{1}{2}}_{2,1}}\,d\tau\\
&\lesssim\big(\sup_{0\leq\tau\leq t}\langle \tau\rangle^{\frac 32}\|c(\tau)\|_{\dot  B^{\frac{3}{2}}_{2,1}}^{\ell}\big)
\int_0^t\langle t-\tau\rangle^{-(\frac 34+\frac s2)}\langle \tau\rangle^{-\frac{3}{2}}\| H\|_{\dot  B^{\frac{3}{2}}_{2,1}}
\|\nabla H\|_{\dot  B^{-\frac{1}{2}}_{2,1}}\,d\tau\\
&\quad+\big(\sup_{0\leq\tau\leq t}\langle \tau\rangle^{\alpha}\|c(\tau)\|_{\dot  B^{\frac{3}{2}}_{2,1}}^{h}\big)
\int_0^t\langle t-\tau\rangle^{-(\frac 34+\frac s2)}\langle \tau\rangle^{-\alpha}\| H\|_{\dot  B^{\frac{3}{2}}_{2,1}}
\|\nabla H\|_{\dot  B^{-\frac{1}{2}}_{2,1}}\,d\tau\\
&\lesssim D^{3}(t)\int_0^t\langle t-\tau\rangle^{-(\frac 34+\frac s2)}\langle \tau\rangle^{-\min(4,\alpha+\frac{5}{2},2\alpha+1,3\alpha)}\,d\tau\\
&\lesssim \langle t\rangle^{-(\frac 34+\frac s2)} D^{3}(t).
\end{align*}
Hence\begin{equation*}\label{f123low}
\int_0^t\langle t-\tau\rangle^{-(\frac{3}{4}+\frac{s}{2})} \big\|f_{3}(\tau)\big\|_{\dot B^{-\frac{3}{2}}_{2,\infty}}^\ell d\tau
\lesssim\langle t\rangle^{-(\frac 34+\frac s2)}
\big(X^2(t)+D^2(t)+D^3(t)\big).
\end{equation*}
Thus, we complete the proof of \eqref{f123low1}. Combining  with \eqref{U} and \eqref{f123low1}, we conclude that for all $t\geq0$ and $s\in(-\frac{3}{2},2],$
\begin{equation}
\label{low}
\langle t\rangle^{\frac 34+\frac s2}
 \|(c,u,H)\|^{\ell}_{\dot  B^{s}_{2,1}}\lesssim D_{0}+X^{2}(t)+D^{2}(t)+D^{3}(t).
\end{equation}
\subsubsection*{Step 2: Decay estimates for the high frequencies of $(\nabla c, u,\nabla H)$}
Now,  the starting point is Inequality \eqref{17} which implies that for $q\geq q_0$ and for some
$c_0=c(q_0)>0,$ we have
\begin{align*}
\frac12\frac d{dt}\alpha_q^2+c_0\alpha_q^2
&\leq\Bigl(\|(f_q,g_q,h_q,\nabla f_q,\nabla h_q)\|_{L^2}+\|R_q(u,c)\|_{L^2}
+\|R_q(u,u)\|_{L^2}\\
&\quad+\|R_q(u,H)\|_{L^2}+\|\wt R_k(u,c)\|_{L^2}+\|\wt R_k(u,H)\|_{L^2}+\|\nabla u\|_{L^\infty}\alpha_q\Bigr)\alpha_q,
\end{align*}
where
\begin{align*}
f_q=\ddq f,\quad g_q=\ddq g, \quad h_q=\ddq h,
\end{align*}
in which
$$R_q(u,b)\eqdefa [u\cdot\nabla,\ddq]b=u\cdot\nabla\ddq b-\ddq(u\cdot\nabla b) \quad\quad for \quad b\in\{c,u,H\},$$
$$\wt R_q^i(u,b)\eqdefa [u\cdot\nabla,\partial_i\ddq]b=u\cdot\nabla\partial_i\ddq b-\partial_i\ddq(u\cdot\nabla b)
 \quad\quad for\quad b\in\{c,H\}.$$
After time integration, we discover that
\begin{align*}
e^{c_0t}\alpha_q(t)
\leq&\alpha_q(0)+\int_0^te^{c_0\tau}\Bigl(\|(f_q,g_q,h_q,\nabla f_q,\nabla h_q)\|_{L^2}+\|R_q(u,c)\|_{L^2}+\|R_q(u,u)\|_{L^2}\\
&\quad+\|R_q(u,H)\|_{L^2}+\|\wt R_k(u,c)\|_{L^2}+\|\wt R_k(u,H)\|_{L^2}+\|\nabla u\|_{L^\infty}\alpha_q\Bigr) d\tau.
\end{align*}
For  $q\geq q_0$, we have  $\alpha_q\approx\|(\nabla\ddq c,\ddq u,\nabla\ddq H)\|_{L^2}$. Then,
\begin{align*}
&\langle t\rangle^\alpha\|(\nabla\ddq c,\ddq u,\nabla\ddq H)(t)\|_{L^2}\lesssim \langle t\rangle^\alpha e^{-c_0t}\|(\nabla\ddq c,\ddq u,\nabla\ddq H)(0)\|_{L^2}\\
&\qquad+\langle t\rangle^\alpha\int_0^te^{c_0(\tau-t)}\Bigl(\|(f_q,g_q,h_q,\nabla f_q,\nabla h_q)\|_{L^2}+\|R_q(u,c)\|_{L^2}
+\|R_q(u,u)\|_{L^2}\\
&\qquad+\|R_q(u,H)\|_{L^2}+\|\wt R_k(u,c)\|_{L^2}+\|\wt R_k(u,H)\|_{L^2}+\|\nabla u\|_{L^\infty}\alpha_q\Bigr) d\tau
\end{align*}
and thus, multiplying both sides by $2^{\frac q2},$ taking the supremum on $[0,T],$
and summing up over $q\geq q_0,$
\begin{equation}\label{S7}
\|\langle t\rangle^\alpha(\nabla c, u,\nabla H)\|_{\wt L^\infty_T(\dot B^{\frac 12}_{2,1})}^{h}
\lesssim\|(\nabla c_{0}, u_{0},\nabla H_{0})\|_{\dot B^{\frac 12}_{2,1}}^h
\!+\sum_{q\geq q_0}\sup_{0\leq t\leq T}\Big(\langle t\rangle^\alpha\!\int_0^t\!e^{c_0(\tau-t)}2^{\frac q2}S_q\,d\tau\Big)
\end{equation}
with $S_q\eqdefa\sum_{i=1}^7 S_q^i$ and
$$\displaylines{
S_q^1\eqdefa \|(f_q,g_q,h_q,\nabla f_q,\nabla h_q)\|_{L^2},\quad
S_q^2\eqdefa\|R_q(u,c)\|_{L^2},\quad
S_q^3\eqdefa\|R_q(u,u)\|_{L^2},\cr
S_q^4\eqdefa\|R_q(u,H)\|_{L^2}, \quad\quad
S_q^5\eqdefa\|\wt R_q(u,c)\|_{L^2}, \quad\quad
S_q^6\eqdefa\|\wt R_q(u,H)\|_{L^2}, \cr
S_q^7\eqdefa\|\nabla u\|_{L^\infty}\|(\ddq\nabla c,\ddq u,\ddq\nabla H)\|_{L^2}.}$$
Bounding the sum, for $0\leq t\leq 2$, and taking advantage of Proposition \ref{Pro:1}, we end up with
\begin{equation}
\begin{split}
\label{s7t02}
&\sum_{q\geq q_0}\sup_{0\leq t\leq 2}\Big(\langle t\rangle^\alpha\!\int_0^t\!e^{c_0(\tau-t)}2^{\frac q2}S_q(\tau)\,d\tau\Big)
\lesssim\int_0^2 \sum_{q\geq q_0}2^{\frac q2}S_q(\tau)\,d\tau\\
&\quad\lesssim\int_0^2\big(\|(f,g,h,\nabla f,\nabla h)\|_{\dot B^{\frac 12}_{2,1}}^{h}+\|\nabla u\|_{\dot B^{\frac 32}_{2,1}}\|(c,u,H,\nabla c,\nabla H)\|_{\dot B^{\frac 12}_{2,1}}\big)d\tau\\
&\quad\lesssim\int_0^2\big(\|(\nabla f,g,\nabla h)\|_{\dot B^{\frac 12}_{2,1}}^{h}+\|\nabla u\|_{\dot B^{\frac 32}_{2,1}}\|(c,u,H,\nabla c,\nabla H)\|_{\dot B^{\frac 12}_{2,1}}\big)d\tau\\
&\quad\eqdefa Q_{1}+Q_{2}.
\end{split}
\end{equation}
From Propositions \ref{p26}-\ref{p27},   we estimate  the terms  $Q_{1}$ and $Q_{2}$  as follows
\begin{align*}
\int_0^2\|\nabla f\|_{\dot B^{\frac 12}_{2,1}}^{h}d\tau
&\lesssim \int_0^2\| f\|_{\dot B^{\frac 32}_{2,1}}^{h}d\tau\\
&\lesssim \int_0^2\|c\div u\|_{\dot B^{\frac 32}_{2,1}}^{h}d\tau\\
&\lesssim \int_0^2\|c\|_{\dot B^{\frac 32}_{2,1}} \|\div u\|_{\dot B^{\frac 32}_{2,1}}d\tau\\
&\lesssim \| c\|_{L_{t}^{\infty}(\dot B^{\frac 12\frac 32}_{2,1})} \| u\|_{L_{t}^{1}(\dot B^{\frac 52}_{2,1})}\\
&\lesssim X^{2}(2),
\end{align*}
\begin{align*}
\int_0^2\|g\|_{\dot B^{\frac 12}_{2,1}}^{h}d\tau
&\lesssim\int_0^2\| L_{1}(c)\cA u\|_{\dot B^{\frac 12}_{2,1}}^{h}d\tau+\int_0^2\| L_{2}(c)\nabla c\|_{\dot B^{\frac 12}_{2,1}}^{h}d\tau\\
&\qquad+\int_0^2\Big\|L_{3}(c)\big(\frac{1}{2}\nabla|H|^{2}-H\cdot\nabla H\big)\Big\|_{\dot B^{\frac 12}_{2,1}}^{h}d\tau
\\&\lesssim\int_0^2\|c\|_{\dot B^{\frac 32}_{2,1}}
\|\cA u\|_{\dot B^{\frac 12}_{2,1}}d\tau+\int_0^2\|c\|_{\dot B^{\frac 32}_{2,1}}
\| \nabla c\|_{\dot B^{\frac 12}_{2,1}}d\tau\\
&\quad+\int_0^2\big(1+\|L_{1}(c)\|_{\dot B^{\frac 32}_{2,1}}\big)
\|H\cdot\nabla H\|_{\dot B^{\frac 12}_{2,1}}d\tau\\
&\lesssim\|c\|_{L_{t}^{\infty}(\dot B^{\frac 12\frac 32}_{2,1})}
\| u\|_{L_{t}^{1}(\dot B^{\frac 52}_{2,1})}+\|c\|_{L_{t}^{2}(\dot B^{\frac 32}_{2,1})}^{2}
\\&\quad+\big(1+\|c\|_{L_{t}^{\infty}(\dot B^{\frac 32}_{2,1})}\big)\|H\|_{L_{t}^{2}(\dot B^{\frac 32}_{2,1})}^{2}
\\&\lesssim\|c\|_{L_{t}^{\infty}(\dot B^{\frac 12\frac 32}_{2,1})}
\| u\|_{L_{t}^{1}(\dot B^{\frac 52}_{2,1})}+\|c\|_{L_{t}^{\infty}(\dot B^{\frac 12\frac 32}_{2,1})}
\| c\|_{L_{t}^{1}(\dot B^{\frac 52\frac 32}_{2,1})}
\\&\quad+\big(1+\|c\|_{L_{t}^{\infty}(\dot B^{\frac 32}_{2,1})}\big)\|H\|_{L_{t}^{1}\dot B^{\frac 52\frac 72}_{2,1}}
\|H\|_{L_{t}^{\infty}(\dot B^{\frac 12\frac 32}_{2,1})}
\\&\lesssim X^{2}(2)+X^{3}(2),
\end{align*}
\begin{align*}
&\int_0^2\|\nabla h\|_{\dot B^{\frac 12}_{2,1}}^{h}d\tau
\lesssim \int_0^2\|h\|_{\dot B^{\frac 32}_{2,1}}^{h}d\tau\\
&\quad\lesssim \int_0^2\|H(\div u)\|_{\dot B^{\frac 32}_{2,1}}^{h}d\tau
+\int_0^2\|H\cdot\nabla u\|_{\dot B^{\frac 32}_{2,1}}^{h}d\tau\\
&\qquad+ \int_0^2\Big\|\nabla\times\big(L_{3}(c)(\nabla\times H)\times H\big)\Big\|_{\dot B^{\frac 32}_{2,1}}^{h}d\tau\\
&\quad\lesssim\int_0^2\|H\cdot\nabla u\|_{\dot B^{\frac 32}_{2,1}}^{h}d\tau+ \int_0^2\Big\| \nabla\times\big(L_{3}(c)(\nabla\times H)\times H\big)\Big\|_{\dot B^{\frac 32}_{2,1}}^{h}d\tau\\
&\quad\lesssim\|H\|_{L_{t}^{\infty}(\dot B^{\frac 12\frac 32}_{2,1})}
\| u\|_{L_{t}^{1}(\dot B^{\frac 52}_{2,1})}+\big(1+\|c\|_{L_{T}^{\infty}(\dot B^{\frac 1 2,\frac 3 2}_{2,1})}\big)
\|\nabla H\|_{L_{T}^{1}(\dot B^{\frac 3 2,\frac 5 2}_{2,1})}
\|H\|_{L_{T}^{\infty}(\dot B^{\frac 3 2,\frac 3 2}_{2,1})}\\
&\quad\lesssim X^{2}(2)+X^{3}(2),
\end{align*}
and
\begin{align*}
&\int_0^2\|\nabla u\|_{\dot B^{\frac 32}_{2,1}}
\|(c,u,H,\nabla c,\nabla H)\|_{\dot B^{\frac 12}_{2,1}}d\tau\\
&\quad\lesssim\|\nabla u\|_{L_{T}^{1}(\dot B^{\frac 32}_{2,1})}
\|(c,u,H,\nabla c,\nabla H)\|_{L_{T}^{\infty}(\dot B^{\frac 12}_{2,1})}\\
&\quad\lesssim\|u\|_{L_{T}^{1}(\dot B^{\frac 52}_{2,1})}
\big(\|(c,H)\|_{L_{T}^{\infty}(\dot B^{\frac 12,\frac 32}_{2,1})}
+\|u\|_{L_{T}^{\infty}(\dot B^{\frac 12}_{2,1})}\big)\\
&\quad\lesssim X^{2}(2).
\end{align*}
Therefore, for the case $t\leq2$,
\begin{equation}\label{l2}
\sum_{q\geq q_0}\sup_{0\leq t\leq 2}\Big(
\langle t\rangle^\alpha\!\int_0^t\!e^{c_0(\tau-t)}2^{\frac{q}{2}}S_q\,d\tau\Big)
\lesssim X^2(2)+X^3(2).
\end{equation}
To bound the supremum on $[2,T],$ we split the integral on $[0,t]$ into  integrals
on  $[0,1]$ and $[1,t],$ respectively.
The $[0,1]$ part of the integral is easy to handle, we have
$$\begin{aligned}
\sum_{q\geq q_0}\sup_{2\leq t\leq T}\Big(
\langle t\rangle^\alpha
\!\int_0^1\!e^{c_0(\tau-t)}2^{\frac{q}{2}}S_q(\tau)\,d\tau\Big)
&\leq  \sum_{q\geq q_0}\sup_{2\leq t\leq T} \Big(\langle t\rangle^\alpha e^{-\frac{c_0}2t} \int_0^1 2^{\frac{q}{2}}S_q\,d\tau\Big)\\
&\lesssim \int_0^1  \sum_{q\geq q_0}  2^{\frac{q}{2}}S_q\,d\tau.
\end{aligned}
$$
Hence
\begin{equation}\label{l1}
\sum_{q\geq q_0}\sup_{2\leq t\leq T}
\Big(\langle t\rangle^\alpha\!\int_0^1\!e^{c_0(\tau-t)}2^{\frac q2}S_q(\tau)\,d\tau\Big)\lesssim X^2(1)+X^3(1).
\end{equation}
Let us finally consider the $[1,t]$ part of the integral for $2\leq t\leq T.$  We shall use repeatedly the following inequalities
\begin{equation}\label{eq:135}
\|\tau\nabla u\|_{\wt L^\infty_t(\dot{B}^{\frac 32}_{2,1})}\lesssim D(t),
\end{equation}
and
\begin{equation}\label{eq:136}
\|\tau\nabla H\|_{\wt L^\infty_t(\dot B^{\frac 52}_{2,1})}\lesssim D(t),
\end{equation}
which are straightforward as regards the high frequencies of $u$ and $H$ and stem from
$$
\|\tau\nabla u\|_{\wt L^\infty_t(\dot B^{\frac 32}_{2,1})}^\ell\lesssim
\|\langle\tau\rangle^{\frac 54}\nabla u\|_{\wt L^\infty_t(\dot B^{\frac 32}_{2,1})}^\ell
\lesssim\|\langle\tau\rangle^{\frac 54} u\|_{L^\infty_t(\dot B^{1}_{2,1})}^\ell\leq D(t),
$$
as well as
$$
\|\tau\nabla H\|_{\wt L^\infty_t(\dot B^{\frac 52}_{2,1})}^\ell\lesssim
\|\langle\tau\rangle^{\frac 54}\nabla H\|_{\wt L^\infty_t(\dot B^{\frac 52}_{2,1})}^\ell
\lesssim\|\langle\tau\rangle^{\frac 54} H\|_{L^\infty_t(\dot B^{1}_{2,1})}^\ell\leq D(t),
$$
for the low frequencies of  $u$ and $H$.

Regarding the contribution of $S_q^1,$  by Lemma \ref{lemma2.13} we first notice that
\begin{align*}\label{eq:decay6}
&\sum_{q\geq q_0}\sup_{2\leq t\leq T}\Big(\langle t\rangle^\alpha\!\int_1^t\!e^{c_0(\tau-t)}2^{\frac q2}S_q^1(\tau)\,d\tau\Big)\\
&\quad=\sum_{q\geq q_0}\sup_{2\leq t\leq T}\Big(
\langle t\rangle^\alpha\!\int_1^t\!e^{c_0(\tau-t)}2^{\frac q2}\|(f_q,g_q,h_q,\nabla f_q,\nabla h_q)\|_{L^2}(\tau)\,d\tau\Big)\\
&\quad\lesssim \|\tau^\alpha(f,g,h,\nabla f,\nabla h)\|_{\wt L^\infty_T(\dot B^{\frac 12}_{2,1})}^h\\
&\quad\lesssim \|\tau^\alpha(\nabla f,g,\nabla h)\|_{\wt L^\infty_T(\dot B^{\frac 12}_{2,1})}^h.
\end{align*}
Now, product laws in tilde spaces ensures that
$$
\|\tau^\alpha\nabla f\|_{\wt L^\infty_T(\dot B^{\frac 12}_{2,1})}^h
\lesssim\|\tau^{\alpha-1} c\|_{\wt L^\infty_T(\dot B^{\frac 32}_{2,1})}
\|\tau\div u\|_{\wt L^\infty_T(\dot B^{\frac 32}_{2,1})}.
$$
The high frequencies of the first term of the r.h.s is obviously bounded by $D(T)$. That is,
\begin{equation}
\label{ch32}
\|\tau^{\alpha-1}c\|^h_{\wt L_T^\infty(\dot B^{\frac 32}_{2,1})}
\lesssim \|\tau^{\alpha}c\|^h_{\wt L_T^\infty(\dot B^{\frac 32}_{2,1})}
\lesssim D(T).
\end{equation}
 As for the low frequencies of the first term of the r.h.s, we notice that
for all small enough $\varepsilon>0,$
\begin{equation}\label{cl32}
\begin{split}
\|\tau^{\alpha-1}c\|^\ell_{\wt L_T^\infty(\dot B^{\frac 32}_{2,1})}
&\lesssim
\|\tau^{\alpha-1}c\|^\ell_{L_T^\infty(\dot B^{\frac 32-2\varepsilon}_{2,1})}\\
&\lesssim
\|\tau^{\alpha-\frac 52+\varepsilon}\tau^{\frac 32-\varepsilon}c\|^\ell_{L_T^\infty(\dot B^{\frac 32-2\varepsilon}_{2,1})}\\
&\lesssim
\|\tau^{\frac 32-\varepsilon}c\|^\ell_{L_T^\infty(\dot B^{\frac 32-2\varepsilon}_{2,1})}\\
&\lesssim
 D(T).
\end{split}
\end{equation}
Combining with \eqref{ch32} and \eqref{cl32}, we obtain
\begin{align}
\label{c32}
\|\tau^{\alpha-1}c\|_{\wt L_T^\infty(\dot B^{\frac 32}_{2,1})}
&\lesssim D(T).
\end{align}
Similarly,
\begin{align}
\label{HH32}\|\tau^{\alpha-1} H\|_{\wt L^\infty_T(\dot B^{\frac 32}_{2,1})}
&\lesssim D(T).
\end{align}
Therefore, using \eqref{eq:135} and  \eqref{c32} we get
$$
\|\tau^\alpha\nabla f\|_{\wt L^\infty_T(\dot B^{\frac 12}_{2,1})}^h\lesssim D^2(T).
$$
Noticing that
 $g= L_{2}(c)\nabla c+L_{1}(c)\cA u
- L_{3}(c)\big(\frac{1}{2}\nabla|H|^{2}-H\cdot\nabla H\big)$,
we shall use repeatedly the following inequality
\begin{equation}
\begin{split}
\label{notc32}
\|c\|_{\wt L^\infty_T(\dot B^{\frac 32}_{2,1})}
&\lesssim\|c\|_{\wt L^\infty_T(\dot B^{\frac 12,\frac 32}_{2,1})}\\
&\lesssim X(T).
\end{split}
\end{equation}
Employing  \eqref{ch32} and \eqref{notc32}, we obtain
\begin{align*}
\|\tau^\alpha L_{2}(c)\nabla c^h\|_{\wt L^\infty_T(\dot B^{\frac 12}_{2,1})}
&\lesssim \|c\|_{\wt L^\infty_T(\dot B^{\frac 32}_{2,1})}
\|\tau^\alpha \nabla c\|_{\wt L^\infty_T(\dot B^{\frac 12}_{2,1})}^h\\
&\lesssim X(T)D(T).
\end{align*}
According to   \eqref{cl32} and remembering the definition of $D(t),$  we have
\begin{align*}
\|\tau^\alpha L_{2}(c)\nabla c^\ell\|_{\wt L^\infty_T(\dot B^{\frac 12}_{2,1})}
&\lesssim \|\tau c\|_{\wt L^\infty_T(\dot B^{\frac 32}_{2,1})}\|\tau^{\alpha-1} \nabla c\|^\ell_{\wt L^\infty_T(\dot B^{\frac 12}_{2,1})}\\
&\lesssim \| \tau c\|_{\wt L^\infty_T(\dot B^{\frac 32}_{2,1})}
D(T)\\
&\lesssim \big(\|\tau c\|_{ L^\infty_T(\dot B^{\frac 32-2\varepsilon}_{2,1})}^{\ell}
+\|\tau c\|_{\wt L^\infty_T(\dot B^{\frac 32}_{2,1})}^{h}\big)
D(T)\\
&\lesssim \big(\|\tau^{\frac{3}{2}-\varepsilon} c\|_{ L^\infty_T(\dot B^{\frac 32-2\varepsilon}_{2,1})}^{\ell}
+\|\langle \tau\rangle^{\alpha} c\|_{\wt L^\infty_T(\dot B^{\frac 32}_{2,1})}^{h}\big)
D(T)\\
&\lesssim  D^2(T).
\end{align*}
Thus,
\begin{align*}
\|\tau^\alpha L_{2}(c)\nabla c \|^{h}_{\wt L^\infty_T(\dot B^{\frac 12}_{2,1})}
\lesssim X^{2}(T)+ D^2(T).
\end{align*}
From \eqref{eq:135} and \eqref{c32}, we also see that
\begin{align*}
\|\tau^\alpha L_{1}(c)\cA u\|^{h}_{\wt L^\infty_T(\dot B^{\frac 12}_{2,1})}
&\lesssim\|\tau\nabla^2u\|_{\wt L^\infty_T(\dot B^{\frac 12}_{2,1})}
\|\tau^{\alpha-1} c\|_{\wt L^\infty_T(\dot B^{\frac 32}_{2,1})}\\
&\lesssim D^2(T).
\end{align*}
Employing Propositions \ref{p26}-\ref{p27} and \eqref{HH32}, for the term $L_{3}(c)\big(\frac{1}{2}\nabla|H|^{2}-H\cdot\nabla H\big)$, we have
\begin{equation}
\begin{split}
\label{L31}
&\Big\| \tau^\alpha L_{3}(c)\big(\frac{1}{2}\nabla|H|^{2}-H\cdot\nabla H\big)\Big\|_{\wt L^\infty_T(\dot B^{\frac 12}_{2,1})}^{h}\\
&\quad\lesssim\|L_{3}(c)\|_{\wt L^\infty_T(\dot B^{\frac 32}_{2,1})}
\| \tau^\alpha H\cdot\nabla H\|_{\wt L^\infty_T(\dot B^{\frac 12}_{2,1})}\\
&\quad\lesssim\big(1+\|L_{1}(c)\|_{\wt L^\infty_T(\dot B^{\frac 32}_{2,1})}\big)
\| \tau^{\alpha-1} H\|_{\wt L^\infty_T(\dot B^{\frac 32}_{2,1})}
\| \tau\nabla H\|_{\wt L^\infty_T(\dot B^{\frac 12}_{2,1})}\\
&\quad\lesssim(1+\|L_{1}(c)\|_{\wt L^\infty_T(\dot B^{\frac 32}_{2,1})})
\| \tau^{\alpha-1} H\|_{\wt L^\infty_T(\dot B^{\frac 32}_{2,1})}^{2}\\
&\quad\lesssim  X^{2}(T)+D^{2}(T)+D^{4}(T).
\end{split}
\end{equation}
Remembering  that $h=-H(\div u)+H\cdot\nabla u-\nabla\times\big(L_{3}(c)(\nabla\times H)\times H\big)$,  we have
\begin{equation}
\label{hh}
\begin{split}
&\| \tau^\alpha \nabla h\|_{\wt L^\infty_T(\dot B^{\frac 12}_{2,1})}^{h}
\lesssim\| \tau^\alpha h\|_{\wt L^\infty_T(\dot B^{\frac 32}_{2,1})}^{h}\\
&\quad\lesssim\Big\| \tau^\alpha \Big(-H(\div u)+H\cdot\nabla u-\nabla\times\big(L_{3}(c)(\nabla\times H)\times H\big)\Big)\Big\|_{\wt L^\infty_T(\dot B^{\frac 32}_{2,1})}^{h}\\
&\quad\lesssim  \big\| \tau^\alpha H\cdot\nabla u\big\|_{\wt L^\infty_T(\dot B^{\frac 32}_{2,1})}+
\Big\| \tau^\alpha \nabla\times\big(L_{3}(c)(\nabla\times H)\times H\big)\Big\|_{\wt L^\infty_T(\dot B^{\frac 32}_{2,1})}.
\end{split}
\end{equation}
It is clear that the first term on the right-side of \eqref{hh} may be bounded by virtue of \eqref{eq:135} and \eqref{HH32}, we have
\begin{align*}
\| \tau^\alpha H\cdot\nabla u\|_{\wt L^\infty_T(\dot B^{\frac 32}_{2,1})}
&\lesssim \| \tau^{\alpha-1} H\|_{\wt L^\infty_T(\dot B^{\frac 32}_{2,1})}
\| \tau\nabla u\|\|_{\wt L^\infty_T(\dot B^{\frac 32}_{2,1})}\lesssim D^{2}(T).
\end{align*}
For the second term on the right-side of \eqref{hh},   by virtue of \eqref{eq:136}, \eqref{HH32}  and \eqref{notc32} we obtain
\begin{align*}
&\Big\| \tau^{\alpha}\nabla\times\big(L_{3}(c)(\nabla\times H)\times H)\big)\Big\|_{\wt L^\infty_T(\dot B^{\frac 32}_{2,1})}\\
&\quad\lesssim \Big\| \tau^{\alpha}L_{3}(c)(\nabla\times H)\times H\Big\|_{\wt L^\infty_T(\dot B^{\frac 52}_{2,1})}\\
&\quad\lesssim \|L_{3}(c)\|_{\wt L^\infty_T(\dot B^{\frac 32}_{2,1})}
\big\| \tau^{\alpha}H\cdot\nabla H \big\|_{\wt L^\infty_T(\dot B^{\frac 52}_{2,1})}\\
&\quad\lesssim \|1+L_{1}(c)\|_{\wt L^\infty_T(\dot B^{\frac 32}_{2,1})}
\| \tau^{\alpha-1}H\|_{\wt L^\infty_T(\dot B^{\frac 32}_{2,1})}\| \tau\nabla H \|_{\wt L^\infty_T(\dot B^{\frac 52}_{2,1})}\\
&\quad\lesssim \big(1+\|c\|_{\wt L^\infty_T(\dot B^{\frac 32}_{2,1})}\big)
\|\tau^{\alpha-1}H\|_{\wt L^\infty_T(\dot B^{\frac 32}_{2,1})}\| \tau\nabla H \|_{\wt L^\infty_T(\dot B^{\frac 52}_{2,1})}\\
&\quad\lesssim X^{2}(T)+D^{2}(T)+D^{4}(T).
\end{align*}
We end up with
\begin{align}\label{6.26}
\sum_{q\geq q_0}\sup_{2\leq t\leq T}\Big(
 \tau^\alpha\!\int_1^t\!e^{c_0(\tau-t)}2^{\frac q2}S_q^1(\tau)\,d\tau\Big)
\lesssim X^{2}(T)+D^{2}(T)+D^{4}(T).
\end{align}
To bound the term with $S_q^2,$ we use the fact that
$$
\int_1^te^{c_0(\tau-t)}\|R_q(u,c)\|_{L^2}\,d\tau\leq
\|R_q(\tau u,\tau^{\alpha-1}a)\|_{L^\infty_t(L^2)}
\int_1^te^{c_0(\tau-t)} \tau^{-\alpha}\,d\tau.
$$
Hence, thanks to Lemma \ref{lemma2.13} and Proposition \ref{Pro:1},
\begin{equation*}
\begin{split}
&\sum_{q\geq q_0}\sup_{2\leq t\leq T}
\Big( \tau^\alpha\!\int_1^t\!e^{c_0(\tau-t)}2^{\frac q2}S_q^2(\tau)\,d\tau\Big)\\
&\quad\lesssim\sum_{q\geq q_0}\sup_{2\leq t\leq T}
\Big( \tau^\alpha\!\int_1^t\!e^{c_0(\tau-t)}2^{\frac q2}\|R_q(u,c)\|_{L^2}\,d\tau\Big)\\
&\quad\lesssim \|\tau\nabla u\|_{\wt L^\infty_T(\dot B^{\frac 32}_{2,1})} \|\tau^{\alpha-1}c\|_{\wt L^\infty_T(\dot B^{\frac 12}_{2,1})}.
\end{split}
\end{equation*}
The first term on the right-side of the above inequality may be bounded thanks to \eqref{eq:135},
and the high frequencies of the last one on the right-side are obviously bounded by $D(T).$
To bound the term   $\|\tau^{\alpha-1}c\|_{\wt L^\infty_T(\dot B^{\frac 12}_{2,1})}^\ell,$
we have the following inequality
\begin{equation*}
\begin{split}
\label{c111}
\|\tau^{\alpha-1}c\|_{\wt L^\infty_T(\dot B^{\frac 12}_{2,1})}^\ell
&\lesssim\|\tau^{\alpha-1}c\|_{L^\infty_T(\dot B^{\frac 12-2\varepsilon}_{2,1})}^\ell\\
&\lesssim \|\tau^{\alpha-2+\varepsilon}\tau^{1-\varepsilon}c\|_{L^\infty_T(\dot B^{\frac 12-2\varepsilon}_{2,1})}^\ell\\
&\lesssim D(T).
\end{split}
\end{equation*}
We eventually get
\begin{equation}
\begin{split}
\label{6.27}
\sum_{q\geq q_0}\sup_{2\leq t\leq T}
\Bigl(\tau^\alpha\!\int_1^t\!e^{c_0(\tau-t)}2^{\frac q2}S_q^2(\tau)\,d\tau\Big)
\lesssim D^2(T).
\end{split}
\end{equation}
Similarly,  we have
\begin{equation}
\label{6.28}
\begin{split}
\sum_{q\geq q_0}\sup_{2\leq t\leq T}
\Big(\tau^\alpha\!\int_1^t\!e^{c_0(\tau-t)}2^{\frac q2}\big(S_q^3(\tau)+S_q^4(\tau)\big)\,d\tau\Big)\lesssim D^2(T).
\end{split}
\end{equation}
Finally, using product laws, \eqref{eq:135}, \eqref{c32}, \eqref{HH32} and Lemma \ref{lemma2.13}, we obtain
\begin{equation}\label{6.29}
\begin{split}
&\sum_{q\geq q_0}\sup_{2\leq t\leq T}
\Big( t^\alpha\!\int_1^t\!e^{c_0(\tau-t)}2^{\frac q2}S_q^5(\tau)\,d\tau\Big)\\
&\quad\lesssim\sum_{q\geq q_0}\sup_{2\leq t\leq T}
\Big( t^\alpha\!\int_1^t\!e^{c_0(\tau-t)}2^{\frac q2}
\|\wt R_q(u,c)\|_{L^2}\,d\tau\Big)\\
&\quad\lesssim \|\tau\nabla u\|_{\wt L^\infty_T(\dot B^{\frac 32}_{2,1})}\|\tau^{\alpha-1}\nabla c \|_{\wt L^\infty_T(\dot B^{\frac 12}_{2,1})}
\sup_{2\leq t\leq T}\Big(  t^\alpha \int_1^te^{c_0(\tau-t)}\tau^{-\alpha}\,d\tau\Big)\\
&\quad\lesssim D^2(T),
\end{split}
\end{equation}
\begin{equation}\label{6.30}
\begin{split}
&\sum_{q\geq q_0}\sup_{2\leq t\leq T}
\Big( t^\alpha\!\int_1^t\!e^{c_0(\tau-t)}2^{\frac q2}S_q^6(\tau)\,d\tau\Big)\\
&\quad\lesssim\sum_{q\geq q_0}\sup_{2\leq t\leq T}
\Big( t^\alpha\!\int_1^t\!e^{c_0(\tau-t)}2^{\frac q2}
\|\wt R_q(u,H)\|_{L^2}\,d\tau\Big)\\
&\quad\lesssim \|\tau\nabla u\|_{\wt L^\infty_T(\dot B^{\frac 32}_{2,1})}\|\tau^{\alpha-1}\nabla H \|_{\wt L^\infty_T(\dot B^{\frac 12}_{2,1})}
\sup_{2\leq t\leq T}\Big( t^\alpha \int_1^te^{c_0(\tau-t)}\tau^{-\alpha}\,d\tau\Big)\\
&\quad\lesssim D^2(T),
\end{split}
\end{equation}
and
\begin{equation}\label{6.31}
\begin{split}
&\sum_{q\geq q_0}\sup_{2\leq t\leq T}\Big( t^\alpha\!\int_1^t\!e^{c_0(\tau-t)}2^{\frac q2}S_k^7(\tau)\,d\tau\Big)\\
&\quad\lesssim\sum_{q\geq q_0}\sup_{2\leq t\leq T}\Big( t^\alpha\!\int_1^t\!e^{c_0(\tau-t)}2^{\frac q2}\|\nabla u\|_{L^\infty}\|(\ddq\nabla c,\ddq u,\ddq\nabla H)\|_{L^2}(\tau)\,d\tau\Big)\\
&\quad\lesssim \|\tau\nabla u\|_{\wt L^\infty_T(\dot B^{\frac 32}_{2,1})}\|\tau^{\alpha-1}(\nabla c,u,\nabla H)\|_{\wt L^\infty_T(\dot B^{\frac 12}_{2,1})}
\sup_{2\leq t\leq T}\Big(  t^\alpha \int_1^te^{c_0(\tau-t)}\tau^{-\alpha}\,d\tau\Big)\\
&\quad\lesssim D^2(T).
\end{split}
\end{equation}
Putting all the above inequalities \eqref{6.26}-\eqref{6.31} together, we conclude that
\begin{equation}\label{6.32}
\sum_{q\geq q_0}\sup_{2\leq t\leq T}\Big(t^\alpha\!\int_1^t\!e^{c_0(\tau-t)}2^{\frac q2}S_q(\tau)\,d\tau\Big)
\lesssim X^{2}(T)+X^{3}(T)+D^{2}(T)+D^{4}(T).
\end{equation}
Then plugging \eqref{l2}, \eqref{l1} and \eqref{6.32} into  \eqref{S7} yields
\begin{equation}\label{high}
\|\langle\tau\rangle^\alpha(\nabla c,u,\nabla H)\|^h_{\wt L^\infty_T(\dot B^{\frac 12}_{2,1})}\lesssim
\|(\nabla c_{0}, u_{0},\nabla H_{0})\|^h_{\dot B^{\frac 12}_{2,1}}+X^2(T)+X^3(T)+D^2(T)+D^4(T).
\end{equation}
\subsubsection*{Step 3: Decay estimates with gain of regularity  for the high frequencies of $\nabla u, \Delta H$.}
This step is devoted to bounding the last  two terms of $D(t)$. We first deal with the term $\|\tau\nabla  u\|_{\wt L^\infty_t(\dot B^{\frac{3}{2}}_{2,1})}^h$  and shall use the fact that the velocity $u$ satisfies
\begin{equation}
\begin{split}
\partial_tu-\cA u=F:=-\nabla c-u\cdot\nabla u-L_{1}(c)\cA u+ L_{2}(c)\nabla c
- L_{3}(c)\big(\frac{1}{2}\nabla|H|^{2}-H\cdot\nabla H\big).
\end{split}
\end{equation}
So,
$$
\partial_t(t\cA u)-\cA(t\cA u)=\cA u+t\cA F.
$$
We deduce from Remark \ref{2.14}  that
$$
\|\tau\cA u\|_{\wt L_t^\infty(\dot B^{\frac 12}_{2,1})}^h
\lesssim \|\cA u\|_{L_t^1(\dot B^{\frac 12}_{2,1})}^h +\|\tau\cA F\|_{\wt L^\infty_t(\dot B^{-\frac 32}_{2,1})}^h,
$$
whence, using the bounds given by Theorem \ref{th:main1} ,
\begin{equation}\label{step3}
\|\tau\nabla u\|_{\wt L_t^\infty(\dot B^{\frac 32}_{2,1})}^h
\lesssim  X(0)+\|\tau F\|_{\wt L^\infty_t(\dot B^{\frac 12}_{2,1})}^h.
\end{equation}
In order to bound the first term of $F$, we notice that, because $\alpha\geq1$ and according to \eqref{high}, we have
$$
\|\tau\nabla c\|_{\wt L^\infty_t(\dot B^{\frac 12}_{2,1})}^h\lesssim
\|\langle\tau\rangle^\alpha c\|_{\wt L^\infty_t(\dot B^{\frac 32}_{2,1})}^h
\lesssim X(0)+X^2(t)+X^3(t)+D^2(t)+D^4(t).$$
Furthermore, from \eqref{eq:135} and   the definition of $X(t)$, we have
\begin{align*}
\|\tau \, u\cdot\nabla u\|_{\wt L^\infty_t(\dot B^{\frac 12}_{2,1})}^h
&\lesssim \|u\|_{\wt L^\infty_t(\dot B^{\frac 12}_{2,1})}
\|\tau \nabla u\|_{\wt L^\infty_t(\dot B^{\frac 32}_{2,1})}\\
&\lesssim X(t)D(t)
\end{align*}
and
\begin{align*}
&\|\tau\, L_{1}(c)\cA u\|_{\wt L^\infty_t(\dot B^{\frac 12}_{2,1})}^h\\
&\lesssim \| c\|_{\wt L^\infty_t(\dot B^{\frac 32}_{2,1})}
 \|\tau\nabla u\|_{\wt L^\infty_t(\dot B^{\frac 32}_{2,1})}\\
&\lesssim X(t)D(t).
\end{align*}
Next, product and composition estimates adapted to tilde spaces give
$$
\|\tau\, L_{2}(c)\nabla c\|_{\wt L^\infty_t(\dot B^{\frac 12}_{2,1})}^h\lesssim \|\tau^{\frac12} c\|_{\wt L^\infty_t(\dot B^{\frac 32}_{2,1})}^2\lesssim D^2(t).
$$
Employing \eqref{HH32} and the definition of $X(t)$, we get
\begin{align*}
&\Big\|\tau\, L_{3}(c)\big(\frac{1}{2}\nabla|H|^{2}-H\cdot\nabla H\big)\Big\|_{\wt L^\infty_t(\dot B^{\frac 12}_{2,1})}^h\\
&\quad\lesssim \| L_{3}(c)\|_{\wt L^\infty_t(\dot B^{\frac 32}_{2,1})}
\|\tau H\cdot\nabla H\|_{\wt L^\infty_t(\dot B^{\frac 12}_{2,1})}\\
&\quad\lesssim \big(1+\| L_{1}(c)\|_{\wt L^\infty_t(\dot B^{\frac 32}_{2,1})}\big)
\| H\|_{\wt L^\infty_t(\dot B^{\frac 32}_{2,1})}
\|\tau \nabla H\|_{\wt L^\infty_t(\dot B^{\frac 12}_{2,1})}\\
&\quad\lesssim \big(1+\|c\|_{\wt L^\infty_t(\dot B^{\frac 32}_{2,1})}\big)
\| H\|_{\wt L^\infty_t(\dot B^{\frac 32}_{2,1})}
\|\tau^{\alpha-1} \nabla H\|_{\wt L^\infty_t(\dot B^{\frac 12}_{2,1})}\\
&\quad\lesssim X^{2}(t)+D^{2}(t)+X^{4}(t).
\end{align*}
Therefore,
\begin{equation}\label{last two}
\|\tau\nabla u\|_{\wt L_t^\infty(\dot B^{\frac 32}_{2,1})}^h\\
\lesssim X(0)+X^2(t)+X^3(t)+X^{4}(t)+D^2(t)+D^4(t).
\end{equation}
Finally, in order to bound the  term of $\|\tau\nabla  H\|_{\wt L^\infty_t(\dot B^{\frac{5}{2}}_{2,1})}^h$,  we shall use the fact that the magnetic $H$ satisfies
\begin{equation}
\begin{split}
\partial_tH-\Delta H=G:=-u\cdot\nabla H-H(\div u)+H\cdot\nabla u-\nabla\times\big(L_{3}(c)(\nabla\times H)\times H\big).
\end{split}
\end{equation}
Furthermore,
$$
\partial_t(t\nabla^{2}H)-\Delta(t\nabla^{2}H)=\nabla^{2}H+t\nabla^{2}G.
$$
From Proposition \ref{Pro:3},   we have
$$
\|\tau\nabla^{2} H\|_{\wt L_t^\infty(\dot B^{\frac 32}_{2,1})}^h
\lesssim \|\nabla^{2} H\|_{L_t^1(\dot B^{\frac 32}_{2,1})}^h +\|\tau\nabla^{2} G\|_{\wt L^\infty_t(\dot B^{-\frac 12}_{2,1})}^h.
$$
Using the bounds given by Theorem \ref{th:main1}, we get
\begin{equation}\label{HH}
\|\tau\nabla^{2} H\|_{\wt L_t^\infty(\dot B^{\frac 32}_{2,1})}^h
\lesssim  X(0) +\|\tau G\|_{\wt L^\infty_t(\dot B^{\frac 32}_{2,1})}^h.
\end{equation}
Bounding the first two terms of $G$, we notice that
\begin{equation}
\begin{split}
\|\tau u\cdot\nabla H\|_{\wt L^\infty_t(\dot B^{\frac 32}_{2,1})}^h
&\lesssim \|u\|_{\wt L^\infty_t(\dot B^{\frac 32}_{2,1})}
\|\tau\nabla H\|_{\wt L^\infty_t(\dot B^{\frac 32}_{2,1})}\\
&\lesssim \|u\|_{\wt L^\infty_t(\dot B^{\frac 12\frac 32}_{2,1})}
\|\tau\nabla H\|_{\wt L^\infty_t(\dot B^{\frac 32}_{2,1})}\\
&\lesssim X(t)\big(\|\tau\nabla H\|_{\wt L^\infty_t(\dot B^{\frac 32}_{2,1})}^{\ell}+
\|\tau\nabla H\|_{\wt L^\infty_t(\dot B^{\frac 32}_{2,1})}^h\big)\\
&\lesssim X(t)\big(\|\tau^{2-\varepsilon} H\|_{ L^\infty_t(\dot B^{\frac 52-2\varepsilon}_{2,1})}^{\ell}+
\|\tau\nabla H\|_{\wt L^\infty_t(\dot B^{\frac 52}_{2,1})}^h\big)\\
&\lesssim X(t)D(t)
\end{split}
\end{equation}
and
\begin{equation}
\begin{split}
\|\tau H\cdot\nabla u\|_{\wt L^\infty_t(\dot B^{\frac 32}_{2,1})}^h
&\lesssim \|H\|_{\wt L^\infty_t(\dot B^{\frac 32}_{2,1})}
\|\tau\nabla u\|_{\wt L^\infty_t(\dot B^{\frac 32}_{2,1})}\\
&\lesssim \|H\|_{\wt L^\infty_t(\dot B^{\frac 12\frac 32}_{2,1})}
\|\tau\nabla u\|_{\wt L^\infty_t(\dot B^{\frac 32}_{2,1})}\\
&\lesssim X(t)D(t).
\end{split}
\end{equation}
The third term of $G$ is similar to the second one, using \eqref{HH32}, we obtain
\begin{align*}
&\Big\|\tau\nabla\times\big(L_{3}(c)(\nabla\times H)\times H)\big)\Big\|_{\wt L^\infty_T(\dot B^{\frac 32}_{2,1})}^{h}\\
&\quad\lesssim \Big\|\tau\nabla\times\big(L_{3}(c)(\nabla\times H)\times H)\big)\Big\|_{\wt L^\infty_T(\dot B^{\frac 12,\frac 32}_{2,1})}\\
&\quad\lesssim \Big\|\tau(L_{3}(c)(\nabla\times H)\times H))\Big\|_{\wt L^\infty_T(\dot B^{\frac 32,\frac 52}_{2,1})}\\
&\quad\lesssim \|L_{3}(c)\|_{\wt L^\infty_T(\dot B^{\frac 32}_{2,1})}
\|\tau H\cdot\nabla H \|_{\wt L^\infty_T(\dot B^{\frac 32,\frac 52}_{2,1})}\\
&\quad\lesssim \|1+L_{1}(c)\|_{\wt L^\infty_T(\dot B^{\frac 32}_{2,1})}
\| H\|_{\wt L^\infty_T(\dot B^{\frac 32}_{2,1})}
\|\tau\nabla H \|_{\wt L^\infty_T(\dot B^{\frac 32,\frac 52}_{2,1})}\\
&\quad\lesssim \big(1+\|c\|_{\wt L^\infty_T(\dot B^{\frac 32}_{2,1})}\big)
\| H\|_{\wt L^\infty_T(\dot B^{\frac 32}_{2,1})}
\|\tau\nabla H \|_{\wt L^\infty_T(\dot B^{\frac 12,\frac 52}_{2,1})}\\
&\quad\lesssim X^{2}(t)+D^{2}(t)+X^{4}(t).
\end{align*}
Hence, reverting to \eqref{HH}, we get
\begin{equation}\label{eq:decay8}
\|\tau\nabla^{2} H\|_{\wt L_t^\infty(\dot B^{\frac 32}_{2,1})}^h\\
\lesssim  X(0)
 +X^{2}(t)+D^{2}(t)+X^{4}(t).
\end{equation}
Finally, adding up the obtained inequality to \eqref{low}, \eqref{high} and \eqref{last two} yields for all $t\geq0,$
\begin{align*}D(t)&\lesssim  X(0)+D_{0}+\|(\nabla c_{0}, u_{0},\nabla H_{0})\|^h_{\dot B^{\frac 12}_{2,1}}+X^{2}(t)+X^{3}(t)+X^{4}(t)+D^{2}(t)+D^{3}(t)+D^{4}(t)\\
&\lesssim  D_{0}+\|(\nabla c_{0}, u_{0},\nabla H_{0})\|^h_{\dot B^{\frac 12}_{2,1}}+X^{2}(t)+X^{3}(t)+X^{4}(t)+D^{2}(t)+D^{3}(t)+D^{4}(t),
\end{align*}
where we have used $X(0)^{\ell}=\|( c_{0}, u_{0}, H_{0})\|^{\ell}_{\dot B^{\frac 12}_{2,1}}\lesssim \|(c_0,u_0,H_0)\|^{\ell}_{\dot B^{-\frac 32}_{2,\infty}}$.
As Theorem \ref{th:main1} ensures that $X(t)$ is small, one can conclude that \eqref{1.2} is fulfilled for all time if $D_{0}$ and $\|(\nabla c_{0}, u_{0},\nabla H_{0})\|^h_{\dot B^{\frac 12}_{2,1}}$ are small enough. This completes the proof of Theorem \ref{th:decay}.

\begin{center}

\end{center}
\end{document}